\documentclass[12pt,twoside]{amsart}
\usepackage{amssymb}
\usepackage{amscd}
\usepackage{xypic}

\title[Injectivity Theorem II]
{A transcendental approach to Koll\'ar's injectivity theorem II}
\author{Osamu Fujino} 
\date{2012/3/5, version 1.35}
\address{Department of Mathematics, Faculty of Science, Kyoto University, 
Kyoto 606-8502, Japan} 
\email{fujino@math.kyoto-u.ac.jp}
\newcommand{\lla}[0]{{\langle\!\langle}}
\newcommand{\rra}[0]{{\rangle\!\rangle}}
\newcommand{\Hom}[0]{{\operatorname{Hom}}}
\newcommand{\loc}[0]{{\operatorname{loc}}}
\newcommand{\Pic}[0]{{\operatorname{Pic}}}
\newcommand{\Supp}[0]{{\operatorname{Supp}}}

\newcommand{\Ker}[0]{{\operatorname{Ker}}}
\newcommand{\xId}[0]{{\operatorname{Id}}}
\newcommand{\Nak}[0]{{\operatorname{Nak}}}
\newcommand{\rank}[0]{{\operatorname{rank}}}
\newcommand{\xre}[0]{{\operatorname{Re}}}
\newcommand{\xIm}[0]{{\operatorname{Im}}}
\newcommand{\Dom}[0]{{\operatorname{Dom}}}
\newtheorem{thm}{Theorem}[section]
\newtheorem{lem}[thm]{Lemma}
\newtheorem{cor}[thm]{Corollary}
\newtheorem{prop}[thm]{Proposition}
\newtheorem*{claim}{Claim}
\newtheorem{problem}[thm]{Problem}

\theoremstyle{definition}
\newtheorem{ques}[thm]{Question}
\newtheorem{ex}[thm]{Example}
\newtheorem{defn}[thm]{Definition}
\newtheorem{rem}[thm]{Remark}
\newtheorem*{ack}{Acknowledgments}       
\newtheorem{say}[thm]{}
\begin{document}
\bibliographystyle{amsalpha+}

\begin{abstract}
We treat a relative version of the main theorem in \cite{f2}:~A 
transcendental approach to Koll\'ar's injectivity theorem. 
More explicitly, we give a curvature condition that implies Koll\'ar 
type cohomology injectivity theorems in the relative setting. 
To carry out this generalization, we use the Ohsawa--Takegoshi 
twisted version of Nakano's identity. 
\end{abstract}

\maketitle
\tableofcontents
\section{Introduction}\label{int}
The following theorem is the main theorem of this paper, 
which is a relative version of the main 
theorem in \cite{f2}. It is a generalization 
of Koll\'ar's injectivity theorem (cf.~\cite[Theorem 2.2]{kollarI}). 
More precisely, it is a generalization of Enoki's injectivity theorem, which is 
an analytic version of Koll\'ar's theorem (see \cite[Theorem 0.2]{enoki} and 
\cite[Corollary 1.4]{f2}).  
We note that Koll\'ar's proofs and the approach of 
Esnault--Viehweg are geometric (cf.~\cite{kollarI}, 
\cite[Chapters 9 and 10]{kollar-book}, and \cite{ev}) 
and are not related to curvature conditions. 
Therefore, we do not know the true relationship between 
Koll\'ar's injectivity theorem and Enoki's one. 
All the terms in the statements will be defined in Section \ref{sec2}. 

\begin{thm}[Main Theorem]\label{II}
Let $f:X\to Y$ be a proper surjective morphism from 
a K\"ahler manifold $X$ to a complex variety $Y$. 
Let $(E, h_E)$ $($resp.~$(L, h_L)$$)$ be a holomorphic vector $($resp.~line$)$ 
bundle on $X$ with a smooth hermitian metric $h_E$ $($resp.~$h_L$$)$. 
Let $F$ be a holomorphic line bundle on $X$ with a singular hermitian 
metric $h_F$. 
Assume the following conditions. 
\begin{itemize}
\item[(i)] There exists a subvariety $Z$ of $X$ such that 
$h_F$ is smooth on $X\setminus Z$. 
\item[(ii)] $\sqrt{-1}\Theta (F)\geq -\widetilde \gamma$ in the sense of 
currents, where $\widetilde\gamma$ is a smooth 
$(1,1)$-form on $X$. 
\item[(iii)] $\sqrt{-1}(\Theta (E)+\xId_E\otimes \Theta
(F))\geq _{\Nak}0$ on $X\setminus Z$. 
\item[(iv)] $\sqrt {-1} (\Theta (E)+\xId_E \otimes \Theta(F)
-\varepsilon_0 \xId _E\otimes \Theta (L))\geq _{\Nak}0$ on 
$X\setminus Z$ for some positive constant $\varepsilon_0$. 
\end{itemize}
Here, 
$\geq _{\Nak}0 $ means the Nakano semi-positivity. 
Let $s$ be a nonzero holomorphic section of $L$. 
Then the multiplication homomorphism 
$$
\times s : 
R^qf_*(K_X\otimes E\otimes F\otimes \mathcal J(h_F))\to 
R^qf_*(K_X\otimes E\otimes F\otimes \mathcal J(h_F)\otimes L)
$$ 
is injective for every $q\geq 0$, where 
$K_X$ is the canonical line bundle of 
$X$ and $\mathcal J(h_F)$ is the 
multiplier ideal sheaf associated to the singular hermitian metric 
$h_F$ of $F$. 
Note that $\times s$ is the sheaf homomorphism induced by 
the tensor product with $s$. 
\end{thm} 

We note that Theorem \ref{II} will be 
generalized slightly in Proposition \ref{326-ne} below. 
For the absolute case and the background of Koll\'ar type cohomology injectivity 
theorems, see the introduction of \cite{f2}. 
The reader who reads Japanese may find \cite{f3} also useful. 
The essential part 
of Theorem \ref{II} is 
contained in Ohsawa's injectivity theorem (see \cite{ohsawa}). 
Our formulation is much more 
suitable for geometric applications than Ohsawa's (cf.~\cite[4.~Applications]{f2}). 
We note that 
the main ingredient of our proof is the Ohsawa--Takegoshi twisted version of 
Nakano's identity (cf.~Proposition \ref{twist}). 

The next corollary directly follows from Theorem \ref{II}. 
It contains a generalization of the Grauert--Riemenschneider 
vanishing theorem. 

\begin{cor}
[Torsion-freeness]\label{to}
Let $f:X\to Y$ be a proper surjective morphism from 
a K\"ahler manifold $X$ to a complex variety $Y$. 
Let $(E, h_E)$ $($resp.~$(F, h_F)$$)$ be a holomorphic vector $($resp.~line$)$ 
bundle on $X$ with a smooth hermitian metric $h_E$ 
$($resp.~a singular hermitian metric $h_F$$)$. 
Assume the conditions {\em{(i)}}, {\em{(ii)}}, and {\em{(iii)}} 
in {\em{Theorem \ref{II}}}. 
Then, 
$R^qf_*(K_X\otimes E\otimes F\otimes \mathcal J(h_F))$ 
is torsion-free for every $q\geq 0$. In particular, 
$R^qf_*(K_X\otimes E\otimes F\otimes \mathcal J(h_F))=0$ 
for $q>\dim X-\dim Y$. 
\end{cor}

We will describe the proof of Theorem \ref{II} in Section \ref{sec3}, 
which may help the reader to understand \cite{ohsawa}. 

From now on, we discuss various vanishing theorems as applications of Theorem 
\ref{II}. 

\begin{cor}[Kawamata--Viehweg--Nadel type vanishing theorem]\label{kvn-vani} 
Let $f:X\to Y$ be a proper surjective morphism from a complex manifold $X$ to a complex 
variety $Y$. 
Let $E$ be a Nakano semi-positive vector bundle on 
$X$ and let $\mathcal L$ be a holomorphic line bundle on $X$ such 
that $\mathcal L^{\otimes m} \simeq \mathcal M\otimes \mathcal O_X(D)$ 
where $m$ is a positive integer, $\mathcal M$ is an $f$-nef-big line bundle, 
and $D$ is 
an effective Cartier divisor on $X$. 
Then $$R^qf_*(K_X\otimes E\otimes \mathcal L\otimes \mathcal J)=0$$ for 
every $q>0$ 
where $\mathcal J=\mathcal J(\frac{1}{m}D)$ is the multiplier ideal 
sheaf associated to $\frac{1}{m}D$.  
\end{cor}
In the minimal model program for projective morphisms 
between complex varieties, the Kawamata--Viehweg vanishing theorem 
plays crucial roles. 
It was first obtained by Nakayama (cf.~\cite[Theorem 3.7]{nakayama1}). 

\begin{cor}[Kawamata--Viehweg vanishing theorem 
for proper morphisms]\label{kv-vani}  
Let $f:X\to Y$ be a proper surjective morphism 
from a complex manifold $X$ to a complex 
variety $Y$. 
Let $H$ be a $\mathbb Q$-Cartier $\mathbb Q$-divisor on $X$ 
such that $\Supp \{H\}$ is a normal crossing divisor on $X$ and that 
$\mathcal O_X(mH)$ is an $f$-nef-big line bundle for some positive integer $m$. 
Then $$R^qf_*(K_X\otimes \mathcal O_X(\ulcorner H\urcorner))=0$$ for 
every $q>0$.  
\end{cor}

We can also prove a Koll\'ar type vanishing theorem from Theorem \ref{II}. 
The proof of Corollary \ref{ko-vani} is a routine work for experts. 

\begin{cor}[Koll\'ar type vanishing theorem]\label{ko-vani} 
Let $f:X\to Y$ be a proper surjective morphism from a 
K\"ahler manifold $X$ to a complex variety $Y$. Let $g:Y\to Z$ 
be a proper morphism between complex varieties. 
Let $E$ be a Nakano semi-positive vector bundle 
on $X$ and let $\mathcal L$ be a holomorphic line 
bundle on $X$ such that $\mathcal L^{\otimes m}\simeq f^*\mathcal N\otimes 
\mathcal O_X(D)$ where $m$ is a positive integer, 
$\mathcal N$ is a $g$-nef-big line bundle, and $D$ is 
an effective Cartier divisor on $X$. Then 
$$
R^pg_*R^qf_*(K_X\otimes E\otimes \mathcal L\otimes \mathcal J)=0
$$ 
for every $p>0$ and every $q\geq 0$ where $\mathcal J=\mathcal J(\frac{1}{m}D)$ 
is the multiplier ideal sheaf associated to $\frac{1}{m}D$. 
\end{cor}

All the statements in the introduction may look complicated. 
So it seems to be worth mentioning that the following well-known 
vanishing theorems easily follow from 
the main theorem:~Theorem \ref{II}. 
The proofs of Corollary \ref{cor-nakano} 
and Corollary \ref{cor-kv} explain one of the reasons 
why we think that injectivity theorems are generalizations 
of vanishing theorems  (cf.~\cite{fujino-multi}, \cite{fujino-vani}) 
and that the formulation of Theorem \ref{II} 
is useful for various applications.  

\begin{cor}[{Nakano vanishing theorem (cf.~\cite[(4.9)]{dem})}]\label{cor-nakano} 
Let $X$ be a compact complex manifold and let $E$ be a 
Nakano positive vector bundle on $X$. 
Then $H^q(X, K_X\otimes E)=0$ for every $q>0$. 
\end{cor}

\begin{proof}
Since $E$ is Nakano positive, $L=\det E$ is a positive line bundle. 
Therefore, $X$ is projective and $L$ is 
ample by Kodaira's embedding theorem. 
If $\varepsilon$ is a small positive number, then 
$$
\sqrt{-1}(\Theta (E)-\varepsilon {\xId}_E\otimes \Theta (L))\geq _{\Nak}0. 
$$
Thus, by using Theorem \ref{II}, $H^q(X, K_X\otimes E)$ can be embedded 
into $H^q(X, K_X\otimes E\otimes L^{\otimes m})$ for 
every sufficiently large positive integer $m$. 
By Serre's vanishing theorem, 
$H^q(X, K_X\otimes E)=0$ for every $q>0$. 
\end{proof}

\begin{cor}[Kawamata--Viehweg vanishing theorem]\label{cor-kv}
Let $X$ be a smooth projective variety and let $L$ be a nef 
and big line bundle on $X$. 
Then $H^q(X, K_X\otimes L)=0$ for every $q>0$. 
\end{cor}

\begin{proof}
By Kodaira's lemma, we can write $L^{\otimes m}
\simeq A\otimes \mathcal O_X(D)$ such that 
$m$ is a positive integer, 
$A$ is an ample line bundle, and 
$D$ is an effective Cartier divisor on $X$ with 
$\mathcal J(\frac{1}{m}D)=\mathcal O_X$.
Let $h_D$ be the singular hermitian metric of $\mathcal O_X(D)$ naturally 
associated to $D$ 
(cf.~\cite[Example 2.3]{f2}) and let $h_A$ be a smooth 
hermitian metric of $A$ whose curvature is positive. 
We put $h_L=h^{\frac{1}{m}}_Ah^{\frac{1}{m}}_D$. 
Then $\sqrt{-1}\Theta (L)\geq 0$ in the sense of currents, $h_L$ is smooth 
on $X\setminus D$, and $\sqrt{-1}(\Theta (L)-\varepsilon \Theta (A))\geq 0$ 
on $X\setminus D$ for 
$0<\varepsilon \ll 1$. 
Therefore, we have inclusions $H^q(X, K_X\otimes L)
\subset H^q(X, K_X\otimes L\otimes A^{\otimes l})$ for every 
$q$ and every sufficiently large positive integer 
$l$ by Theorem \ref{II}. Thus 
$H^q(X, K_X\otimes L)=0$ for every $q>0$ by Serre's vanishing theorem.  
\end{proof}

For related topics, 
see \cite[\S3]{nakayama1}, \cite[II.~\S5.c, V.~\S3]{nakayama2}, and \cite{take}. 
We think that 
one of the most important open problems on vanishing theorems for the minimal model 
program is to prove the results in \cite[Sections 6 and 8]{fujino-fundamental} and 
\cite[Chapter 2]{fujino-book} for projective morphisms between complex 
varieties. The following problem seems to be a first step. 

\begin{problem}\label{prob18} 
Let $f:X\to Y$ be a projective morphism 
from a complex manifold $X$ onto a complex variety $Y$. 
Let $D$ be a simple normal crossing divisor on $X$. 
Let $L$ be a $\pi$-semi-ample Cartier divisor on $X$. 
Let $s$ be a nonzero holomorphic section of $kL$ for some positive integer $k$. 
Assume that 
$(s=0)$ contains no lc centers of $(X, D)$. 
Then 
$$\times s: R^q\pi_*(K_X\otimes \mathcal O_X(D+mL))\to 
R^q\pi_*(K_X\otimes \mathcal O_X(D+(m+k)L)) 
$$
is injective for every $q$ and every positive integer $m$. 
Note that $W\subset X$ is an lc center  of $(X, D)$ if and only if 
$W$ is an irreducible component of $D_{i_1}\cap 
\cdots 
\cap D_{i_l}$ for some 
$\{i_1, \cdots, i_l\}\subset I$, 
where $D=\sum _{i\in I}D_i$ is the irreducible decomposition. 
\end{problem} 

The statement of Problem \ref{prob18} holds true when 
$Y$ is a point (see, for example, \cite{fujino-on}, \cite{fujino-fundamental}, 
\cite{fujino-book}, and \cite{fujino-vani2}). 
We do not know if the statement of Problem \ref{prob18} 
holds true or not under the weaker assumption that 
$\mathcal O_X(L)$ is semi-positive even when 
$Y$ is a point. 
When $D=0$, the statement of Problem \ref{prob18} 
is a special case of Theorem \ref{II}. 
In Problem \ref{prob18}, it may be sufficient to 
assume that $f$ is a K\"ahler morphism. 
For the solution of 
Problem \ref{prob18}, we can not directly use the arguments 
in this paper because the $L^2$-method 
does not work for log canonical pairs. 
We note that the arguments in \cite{fujino-on}, 
\cite{fujino-fundamental}, \cite{fujino-book}, 
\cite{fujino-vani2}, and \cite{fujino-slc}  
are geometric. 

We summarize the contents of this paper. 
In Section \ref{sec2}, we collect basic definitions and 
results in the algebraic and analytic geometries. 
In this section, we discuss the Ohsawa--Takegoshi twisted version of 
Nakano's identity. 
It is a key ingredient of the proof of the main theorem:~Theorem \ref{II}. 
Section \ref{sec3} is devoted to the proof of 
the main theorem:~Theorem \ref{II}. 
In Section \ref{sec-p}, we discuss the proofs 
of the corollaries in Section \ref{int} and 
some applications. 
In the final section:~Section \ref{app}, we 
discuss various examples of nef, semi-positive, and 
semi-ample 
line bundles. 
It is very important to understand the differences in the notion of 
semi-ample, semi-positive, and nef line bundles. 

\begin{ack} 
The first version of this paper was written 
in Nagoya in 2006.  
The author was partially supported by The Sumitomo Foundation 
and by the Grant-in-Aid for Young Scientists (A) 
$\sharp$17684001 from 
JSPS when he prepared the first version. 
He thanks Professor Takeo Ohsawa for giving him 
a preliminary version of \cite{ohsawa}. 
He revised this paper in Kyoto in 2011. 
He was partially supported by The Inamori Foundation and by the 
Grant-in-Aid for Young Scientists (A) $\sharp$20684001 from JSPS. 
Finally, he thanks the referees for many useful comments and 
Shinnosuke Okawa for pointing out a mistake. 
\end{ack} 

\section{Preliminaries}\label{sec2} 

In this section, we collect basic definitions and results in the algebraic 
and analytic geometries. 

\begin{say}[Projective morphisms]
For details of {\em{projective morphisms}} of complex varieties and 
ample line bundles, 
see, for example, \cite[\S1]{nakayama1} and 
\cite[II.~1.10.~Definition, Remark]{nakayama2}. 
\end{say}

\begin{say}[Big line bundles]
In this paper, we will freely use {\em{Iitaka's $D$-dimension}} $\kappa$, 
the {\em{numerical $D$-dimension}} $\nu$, and so on, 
for algebraic varieties. 

Let us recall the definition of {\em{$f$-nef-big line bundles}} for 
proper morphisms between complex varieties, 
which we need in corollaries in Section \ref{int}.

\begin{defn}[{cf.~\cite[Definition]{nakayama1}}] 
Let $f:X\to Y$ be a proper surjective morphism 
from a complex variety $X$ onto a complex variety $Y$. 
Let $L$ be a line bundle on $X$. 
Then $L$ is called {\em{$f$-big}} if the relative Iitaka $D$-dimension 
$\kappa (X/Y, L)=\dim X-\dim Y$. 
Furthermore, if $L\cdot C\geq 0$ for every irreducible 
curve $C$ such that 
$f(C)$ is a point, then $L$ is called {\em{$f$-nef-big}}. 
\end{defn}
\end{say}

\begin{say}[$\mathbb Q$-divisors]
Let $D=\sum _i d_i D_i$ be 
a $\mathbb Q$-divisor on a normal complex variety $X$ where 
$D_i$ is a prime divisor for every $i$ and $D_i\ne D_j$ for 
$i\ne j$. 
Then we define the {\em{round-up}} 
$\ulcorner D\urcorner =\sum _i \ulcorner d_i \urcorner D_i$ 
(resp.~the {\em{round-down}} $\llcorner D\lrcorner 
=\sum _i \llcorner d_i \lrcorner D_i$), 
where for every real number $x$, $\ulcorner x\urcorner$ 
(resp.~$\llcorner x\lrcorner$) 
is the integer defined by $x\leq \ulcorner x\urcorner < x+1$ 
(resp.~$x-1 < \llcorner x\lrcorner \leq x$). 
The {\em{fractional part}} $\{D\}$ 
of $D$ denotes $D-\llcorner D\lrcorner$. 
\end{say}

\begin{say}[Singular hermitian metric] 

Let $L$ be a holomorphic line bundle on a complex manifold $X$. 

\begin{defn}[Singular hermitian metric] 
A {\em{singular hermitian metric}} on $L$ is a metric which is given 
in any trivialization $\theta:L|_{\Omega}\simeq \Omega \times 
\mathbb C$ by 
$$
\| \xi\|=|\theta (\xi)|e^{-\varphi(x)},\ \ \  x\in \Omega, \ \xi \in L_x, 
$$
where $\varphi\in L^1_{\loc}(\Omega)$ is an arbitrary function, 
called the {\em{weight}} of the metric with respect to the trivialization $\theta$. 
Here, $L^1_{\loc}(\Omega)$ is the space of the locally integrable functions 
on $\Omega$. 
\end{defn}
\end{say}

\begin{say}[Multiplier ideal sheaf] 
The notion of multiplier ideal sheaves introduced by Nadel 
is very important. First, we recall the notion 
of (quasi-)plurisubharmonic 
functions. 

\begin{defn}[Plurisubharmonic function] 
A function $u:\Omega\to [-\infty, \infty)$ defined on an open set 
$\Omega\subset \mathbb C^n$ is called {\em{plurisubharmonic}} 
(psh, for short) if 
\begin{itemize}
\item[1.] $u$ is upper semi-continuous, and 
\item[2.] for every complex line $L\subset \mathbb C^n$, 
$u|_{\Omega\cap L}$ is subharmonic on $\Omega\cap L$, that is, 
for every $a\in \Omega$ and $\xi \in \mathbb C^n$ 
satisfying $|\xi|<d(a, \Omega^c)$, the function $u$ satisfies 
the mean inequality 
$$
u(a)\leq \frac{1}{2\pi} \int ^{2\pi}_{0} u(a+e^{i\theta}\xi)d\theta. 
$$
\end{itemize} 
Let $X$ be an $n$-dimensional complex manifold. 
A function $\varphi:X\to [-\infty, \infty)$ is said to be 
{\em{plurisubharmonic}} (psh, for short) if 
there exists an open cover $X=\bigcup _{i\in I}U_i$ such that 
$\varphi|_{U_i}$ is plurisubharmonic on $U_i$ 
($\subset \mathbb C^n$) for every $i$. 
A {\em{smooth strictly plurisubharmonic function}} $\psi$ on $X$ is a 
smooth function on $X$ such that $\sqrt{-1} \partial \bar\partial \psi$ is 
a positive definite smooth $(1,1)$-form. 
\end{defn} 

\begin{defn}
A {\em{quasi-plurisubharmonic}} 
(quasi-psh, for short) function is a function $\varphi$ which 
is locally equal to the sum of a psh function and of a smooth function. 
\end{defn}

Next, we define multiplier ideal sheaves. 

\begin{defn}[Multiplier ideal sheaf]
If $\varphi$ is a quasi-psh function on a complex manifold $X$, 
the {\em{multiplier ideal sheaf}} 
$\mathcal J(\varphi)\subset \mathcal O_X$ is 
defined by 
$$
\Gamma (U, \mathcal J(\varphi))=\{f\in \mathcal O_X(U);\ |f|^2e^{-2\varphi}\in 
L^1_{\loc}(U)\} 
$$ 
for every open set $U\subset X$. 
Then it is known that 
$\mathcal J(\varphi)$ is a coherent ideal sheaf of $\mathcal O_X$. 
See, for example, \cite[(5.7) Proposition]{dem}. 
\end{defn}

Finally, we note the definition of $\mathcal J(h_F)$ in Theorem \ref{II}. 

\begin{rem}
In Theorem \ref{II}, the curvature $\sqrt{-1}\Theta(F)$ of 
$(F, h_F)$ is expressed by 
$\sqrt{-1}\bar\partial \partial\log h_F$. 
Therefore, if we write $h_F=e^{-2\varphi}$ locally, 
then $$\sqrt{-1}\Theta(F)=\sqrt{-1}\bar\partial \partial\log h_F=
2\partial \bar\partial \varphi.$$ 
By the assumption (ii) in Theorem \ref{II}, 
we may assume that the weight $\varphi$ of the singular hermitian metric $h_F$ is a quasi-psh function 
on any trivialization. 
So, we can define multiplier ideal sheaves locally and 
check that they are independent of trivializations. 
Thus, we can define 
the multiplier ideal sheaf globally and denote it by $\mathcal J(h_F)$, 
which is an abuse of notation.  
It is a coherent ideal 
sheaf on $X$. 
\end{rem}
\end{say}

\begin{say}[K\"ahler geometry]  
We collects the basic notion and results 
of the hermitian and K\"ahler geometries (see also \cite{dem}). 

\begin{defn}[Chern connection 
and its curvature form]\label{conecone}
Let $X$ be a complex hermitian manifold and 
let $(E,h)$ be a holomorphic hermitian vector bundle on $X$. 
Then there exists the {\em{Chern connection}} $D=D_{(E, h)}$, 
which can be split in a unique way as a sum of a $(1,0)$ and of 
a $(0,1)$-connection, $D=D'_{(E, h)}+D''_{(E,h)}$. 
By the definition of the Chern connection, 
$D''=D''_{(E,h)}=\bar\partial$. 
We obtain the {\em{curvature form}} $\Theta(E)=\Theta_{(E, h)}
=\Theta_h:=D^2_{(E, h)}$. The subscripts might be suppressed 
if there is no danger of confusion. 

\begin{defn}[Inner product]\label{d216}
Let $X$ be an $n$-dimensional complex 
manifold with the hermitian metric $g$. 
We denote by $\omega$ the {\em{fundamental form}} of $g$. 
Let $(E,h)$ be a hermitian vector bundle on $X$ 
and let $u,v$ be $E$-valued $(p,q)$-forms with measurable coefficients.  
We set 
$$
\|u\|^2=\int _X |u|^2dV_\omega, 
\ \lla u, v \rra= \int _X \langle u, v \rangle 
dV_{\omega}, 
$$ 
where $|u|$ is the pointwise norm induced by $g$ and 
$h$ on $\Lambda^{p,q}T^*_X\otimes E$, and 
$dV_\omega=\frac{1}{n!}\omega^n$. 
More explicitly, $\langle u, v\rangle dV_\omega
={}^t u\wedge H\overline{\ast v}$, where ${}^t u$ 
is the transposed matrix of $u$, $\ast$ 
is the {\em{Hodge star operator}} relative 
to $\omega$, and $H$ is the (local) matrix representation of $h$. 
When we need to emphasize the metrics, we write 
$|u|_{g,h}$, and so on. 
\end{defn}

\begin{lem}[Adjoint]\label{lem-adj}  
Let $\theta\in C^{s,t}(X)$ where 
$C^{s,t}(X)$ is the space of 
smooth $(s, t)$-forms on $X$. 
Let $\theta^*$ be the adjoint operator of 
$\theta\wedge\cdot$ with respect to $\langle\ ,\  \rangle$, that 
is, 
$$
\langle \theta\wedge v, u\rangle=\langle v, \theta^*u\rangle
$$ 
for every $u\in C^{p, q}(X, E)$ and 
every $v\in C^{p-s. q-t}(X, E)$. 
Note that $C^{p, q}(X, E)$ {\em{(}}resp.~$C^{p-s, q-t}(X, E)${\em{)}} is the 
space of smooth $E$-valued $(p, q)$-forms 
{\em{(}}resp.~$(p-s, q-t)$-forms{\em{)}} on $X$. 
Then $\theta^*=(-1)^{(p+q)(s+t+1)}
\ast\bar\theta\ast$ on $C^{p,q}(X, E)$. 
In particular, if $\theta$ 
is a $1$-form, then $\theta^*=\ast \bar\theta \ast$. 
\begin{proof}
We take $u\in C^{p,q}(X, E)$ and $v\in C^{p-s, q-t}(X, E)$. 
Then $$\langle v, \theta^* u\rangle dV_{\omega}=
\langle \theta\wedge v, u\rangle dV_{\omega}$$ by the definition 
of $\theta^*$. 
Let $H$ be the local matrix representation of $h$. 
Then we 
have 
\begin{align*}
\langle \theta\wedge v, u\rangle dV_{\omega}&={}^t(\theta\wedge v)H
\overline{\ast u}\\
&=(-1)^{(p+q-s-t)(s+t)}({}^tv)H\theta \overline {\ast u}\\ 
&=(-1)^{(p+q-s-t)(s+t)+(2n-p-q+s+t)}({}^tv)H \ast \ast 
\overline{\bar\theta \ast u}\\
&=(-1)^{(p+q)(s+t+1)}({}^tv)H\overline {\ast \ast \bar\theta \ast 
u}. 
\end{align*}
Therefore, $\theta ^* u=
(-1)^{(p+q)(s+t+1)}\ast \bar\theta\ast u$. 
\end{proof}
\end{lem}

Let $L^{p,q}_{(2)}(X,E)(=
L^{p.q}_{(2)}(X,(E,h)))$ 
be the space of square integrable $E$-valued $(p,q)$-forms on $X$. 
The inner product was defined in 
Definition \ref{d216}. 
When we emphasize the metrics, 
we write $L^{p,q}_{(2)}(X,E)_{g,h}$, where 
$g$ (resp.~$h$) is the hermitian metric of $X$ 
(resp.~$E$). 
As usual one can view $D'$ and $D''$ as closed and densely 
defined operators on the Hilbert space $L^{p,q}_{(2)}(X,E)$. 
The formal adjoints 
${D'}^*$, ${D''}^*$ also have closed extensions in 
the sense of distributions, which do not necessarily 
coincide with the Hilbert space adjoints in the sense of 
Von Neumann, since 
the latter ones may have strictly smaller domains. 
It is well known, however, that the domains coincide if 
the hermitian metric of $X$ is complete. See 
Lemma \ref{25} below. 

\begin{lem}[Density Lemma]\label{25}
Let $X$ be a complex manifold with the complete hermitian 
metric $g$ and let $(E,h)$ be a holomorphic hermitian vector bundle 
on $X$. 
Then $C^{p,q}_{0}
(X,E)$ is dense in $\Dom (\bar\partial)
\cap \Dom (D''^*_{(E,h)})$ with respect to the 
graph norm $\|v\|+
\|\bar\partial v\|+\|D''^*_{(E,h)}v\|$, 
where $\Dom (\bar\partial)$ $($resp.~$\Dom (D''^*_{(E,h)})$$)$ 
is the domain of $\bar\partial$ $($resp.~$D''^*_{(E,h)}$$)$. 
Here $C^{p, q}_{0}(X, E)$ is the space of smooth 
$E$-valued $(p, q)$-forms with compact supports on $X$.     
\end{lem}
 
Suppose that $(E,h)$ is a 
holomorphic hermitian vector bundle and that $(e_\lambda)$ 
is a holomorphic frame for $E$ over some 
open set $U$. Then the metric $h$ is given by the $r\times r$ 
hermitian matrix $H=(h_{\lambda\mu})$, where $h_{\lambda\mu}=
h(e_\lambda, e_\mu)$ and 
$r=\rank E$. 
Then we have $h(u, v)={}^tuH\bar v$ on $U$ for 
smooth sections $u$, $v$ of $E|_{U}$. 
This implies that $h(u,v)=\sum_{\lambda, \mu}u_\lambda h_{\lambda\mu}
\bar{v}_{\mu}$ for $u=\sum e_iu_i$ and $v=\sum e_j v_j$. 
Then we obtain that 
$\sqrt{-1}\Theta_h(E)
=\sqrt{-1}\bar\partial(\overline{H}^{-1}\partial \overline{H})$ 
and ${}^t\overline{(\sqrt{-1}{}^t\Theta_h(E)H)}=\sqrt{-1}{}^t\Theta_h(E)H$ on 
$U$. 
Let $C^{p,q}(X,E)$ (resp.~$C^{p,q}_0(X,E)$) be the space of smooth 
$E$-valued $(p,q)$-forms (resp.~smooth 
$E$-valued $(p,q)$-forms with compact supports) on $X$. 
We define $\{u,v\}={}^t u\wedge H\bar v$ for $u\in C^{p,q}(X,E)$ 
and $v\in C^{r,s}(X,E)$, where ${}^t u$ is the transposed 
matrix of $u$. We will use $\{\cdot, \cdot\}$ in Section \ref{sec3}. 
\end{defn}

\begin{defn}[Nakano positivity and semi-positivity] 
Let $(E, h)$ be a holomorphic vector bundle on a complex manifold $X$ with 
a smooth hermitian metric $h$. 
Let $\Xi$ be a $\Hom (E, E)$-valued $(1,1)$-form such that 
${}^t\overline{({}^t\Xi h)}={}^t\Xi h$. 
Then $\Xi$ is said to be {\em{Nakano positive}} (resp.~{\em{Nakano 
semi-positive}}) 
if the hermitian form on $T_X\otimes E$ associated to 
${}^t\Xi h$ is positive definite 
(resp.~semi-definite). 
We write $\Xi>_{\Nak}0$ (resp.~$\geq_{\Nak}0$). We note that 
$\Xi_1>_{\Nak}\Xi_2$ (resp.~$\Xi_1\geq _{\Nak}\Xi_2$) means that 
$\Xi_1-\Xi_2>_{\Nak} 0$ (resp.~$\geq _{\Nak} 0$). 
A holomorphic vector bundle $(E, h)$ is said to be {\em{Nakano positive}} 
(resp.~{\em{semi-positive}}) if $\sqrt{-1}\Theta(E)>_{\Nak}0$ (resp.~$\geq _{\Nak}0$). 
We usually omit \lq\lq Nakano\rq\rq when $E$ is a line bundle. 
We often simply say that 
a holomorphic line bundle $L$ is {\em{semi-positive}} 
if there exists a smooth hermitian metric $h_L$ on $L$ such that 
$\sqrt{-1}\Theta (L)\geq 0$. 
\end{defn}

The space of harmonic forms will play important roles in the proof of 
Theorem \ref{II}. See also the introduction of \cite{f2}. 

\begin{defn}[Harmonic forms]\label{harmoni}
Let $X$ be an $n$-dimensional complete K\"ahler manifold 
with a complete K\"ahler metric $g$. 
Let $(E, h_E)$ be a 
holomorphic hermitian vector bundle on $X$. 
We put 
$$
\mathcal H^{p,q}(X, (E, h_E))_g
=\{u\in L^{p, q}_{(2)}(X, E) | \bar\partial u=0 \ {\text{and}} 
\ D''^*_{(E, h_E)}u=0\}. 
$$ 
Note that $
\mathcal H^{p,q}(X, (E, h_E))_g\subset C^{p,q}(X,E)$ by the 
regularization theorem for elliptic partial differential equations of 
second order. 
\end{defn}
\end{say}

\begin{say}[Ohsawa--Takegoshi twist]
The following formula is a {\em{twisted}} version 
of Nakano's identity, which is now 
well known to the experts. 
\begin{prop}[Ohsawa--Takegoshi twist]\label{twist} 
Let $(E, h)$ be a holomorphic hermitian vector 
bundle on an $n$-dimensional K\"ahler manifold $X$. 
Let $\eta$ be any smooth positive function on $X$. 
Then, for every $u\in C^{n,q}_0(X,E)$, 
the equality 
\begin{align*} 
\tag{$\spadesuit$}&\ \ \ \ \ \ \ \ \ \|\sqrt{\eta}D''^*_{(E, h)}u\|^2+
\|\sqrt{\eta}\bar\partial u\|^2-\|\sqrt{\eta}{D'}^*u \|^2\\ 
=&\lla\sqrt{-1}
(\eta \Theta_h-\xId_E\otimes \partial \bar\partial \eta)
\Lambda u, u\rra+2\xre\lla
\bar\partial\eta\wedge D''^*_{(E,h)}u, u\rra
\end{align*}
holds true. 
Here, we denote by $\Lambda$ the adjoint operator 
of $\omega\wedge\cdot\ $. 
Note that $D''=D''_{(E,h)}=\bar\partial$ and 
${D'}^*$ are independent of the hermitian metric $h$. 
\end{prop} 
\begin{proof}[Sketch of the proof] 
We quickly review the proof of this proposition for the reader's convenience. 
If $A, B$ are the endomorphisms of pure degree of the graded module 
$C^{\bullet, \bullet}(X, E)$, their 
{\em{graded Lie bracket}} is defined by 
$$[A, B]=AB-(-1)^{\deg A\deg B}BA.$$ 
Let 
$$\Delta'=D'D'^*+D'^*D'$$ and 
$$\Delta''=D''D''^*+D''^*D''$$ be the complex Laplace operators. 
Then it is well known that 
$$
\Delta''=\Delta'+[\sqrt{-1}\Theta (E), \Lambda],  
$$ 
which is sometimes called Nakano's identity. 
Let us consider the {\em{twisted}} Laplace operators 
$$
D'\eta D'^*+D'^*\eta D'=\eta \Delta'+(\partial \eta)D'^*-(\partial \eta)^*D', 
$$
and 
$$
D''\eta D''^*+D''^*\eta D''
=\eta \Delta''+(\bar \partial \eta )D''^*-(\bar \partial \eta)^*D''. 
$$ 
On the other hand, we can easily check that 
$$
[\sqrt{-1}\partial \bar\partial \eta, \Lambda]=
[D'', (\bar\partial \eta)^*]+[D'^*, \partial \eta] 
$$ 
by $(\bar\partial \eta)^*=-\sqrt{-1}[\partial \eta, \Lambda]$ and 
$D'^*=-\sqrt{-1}[\bar\partial, \Lambda]$. 
Combining these equalities, we find 
\begin{eqnarray*}
&D''\eta D''^*+D''^*\eta D''-D'\eta D'^*-D'^*\eta D'
+[\sqrt{-1}\partial \bar\partial \eta, \Lambda]\\ 
&=\eta[\sqrt{-1}\Theta (E), \Lambda]+(\bar \partial \eta )D''^*
+D''(\bar\partial \eta)^*+(\partial \eta)^*D'+D'^*(\partial \eta). 
\end{eqnarray*}
Apply this identity to a form $u\in C^{n,q}_0(X,E)$ and take the inner 
product with $u$. Then we obtain the desired formula. 
\end{proof}

The next proposition is \cite[Lemma 2.1]{ohsawa}. 
The proof is a routine work. 
It easily follows from Lemmas \ref{25}, 
\ref{26}, and \ref{26.5}. 
\begin{prop}\label{24} 
Fix a complete K\"ahler metric $g$ on $X$. 
We put 
$$
D^{n, q}=\{u\in L^{n,q}_{(2)}(X,E)\ |\ 
\bar\partial u\in L^{n,q+1}_{(2)}(X,E)\ \text{and}\ 
{D''}^*u\in L^{n,q-1}_{(2)}(X,E)\},  
$$ 
that is, $D^{n,q}=\Dom (\bar\partial)
\cap \Dom (D''^*_{(E,h)})\subset L^{n,q}_{(2)}(X,E)$.  
Suppose that $\eta$ is bounded and that there exists 
a constant $\varepsilon>0$ such that 
$$
\sqrt{-1}(\eta\Theta_h-\xId_E\otimes \partial \bar\partial \eta-
\varepsilon \xId_E\otimes \partial 
\eta\wedge \bar\partial \eta)\geq_{\Nak}0
$$ 
holds everywhere. 
Then the equality $(\spadesuit)$ in {\em{Proposition \ref{twist}}} holds for all 
$u\in D^{n,q}$. 
\end{prop}

\begin{lem}\label{26}
For every $u\in C^{n,q}(X,E)$ and 
any positive real number $\delta$, we have 
\begin{eqnarray*}
|2\xre\lla
\bar\partial\eta\wedge D''^*_{(E,h)}u, u\rra|
&=&|2\xre\lla 
D''^*_{(E,h)}u, (\bar\partial \eta)^*u\rra|
\\&\leq& \frac{1}{\delta}\| 
D''^*_{(E,h)}u\|^2+\delta \|(\bar\partial \eta)^*u\|^2, 
\ \text{and}
\end{eqnarray*}
\begin{eqnarray*}
\|(\bar\partial \eta)^*u\|^2&=&\lla 
(\bar\partial \eta)^*u, (\bar\partial \eta)^*u\rra\\ 
&=&\lla \sqrt{-1} \partial \eta\wedge \bar\partial 
\eta\Lambda u, u\rra  
\end{eqnarray*}
since $(\bar\partial \eta)^*u=-\sqrt{-1}\partial \eta \Lambda u$ 
for 
$u\in C^{n,q}(X,E)$. 
Note that $(\bar\partial \eta)^*$ is the adjoint operator of 
$\bar\partial \eta\wedge\cdot$ relative to the inner 
product $\langle\ , \ \rangle$. 
\end{lem}

By combining Proposition \ref{twist} with Lemma \ref{26}, 
we obtain the next lemma. 

\begin{lem}\label{26.5}
We use the same notation as in {\em{Proposition \ref{twist}}}. 
Assume that 
$$
\sqrt{-1}(\eta\Theta_h-\xId_E\otimes \partial \bar\partial \eta-
\varepsilon \xId_E\otimes \partial 
\eta\wedge \bar\partial \eta)\geq_{\Nak}0
$$ 
holds everywhere for some positive constant $\varepsilon$. 
Then, for every $u\in C^{n, q}_0(X, E)$, we have 
\begin{align*}
\|\sqrt{\eta}D'^*u\|^2\leq 
\|\sqrt{\eta}D''^*_{(E, h)}u\|^2+
\|\sqrt{\eta}\bar\partial u\|^2
+\frac{1}{\varepsilon}\|\sqrt{\eta}D''^*_{(E, h)}u\|^2, 
\end{align*}
and 
\begin{align*}
\|\sqrt{\eta}D''^*_{(E, h)}u\|^2+
\|\sqrt{\eta}\bar\partial u\|^2-\|\sqrt{\eta}{D'}^*u \|^2
+\frac{1}{\delta}\|\sqrt{\eta}D''^*_{(E, h)}u\|^2
\\ 
\geq 
(\varepsilon -\delta) \|(\bar\partial \eta)^*u\|^2
\end{align*}
for any positive real number $\delta$. 
\end{lem}

We close this section by the following remark on \cite{take}. 
\begin{rem}
By Proposition \ref{24}, we can prove \cite[Theorem 3.4 (ii)]{take} under 
the slightly weaker assumption that $\varphi$ is a {\em{bounded}} 
smooth psh function on $M$. We do not have to assume that $|d\varphi|$ is 
bounded on $M$. For the notations, see \cite{take}. In this case, there 
are positive constants $C_1$ and $C_2$ such that 
$\varphi+C_1>0$ on $M$ and $C_2-(\varphi+C_1)^2>0$ on $M$. 
We can use Proposition \ref{24} (and Lemmas \ref{26}, \ref{26.5}) for 
$\eta:=C_2-(\varphi+C_1)^2$ and $\varepsilon:=\frac{1}{2C_2}>0$. 
Then we obtain $(\bar\partial \varphi)^*u=0$ and 
$\langle \sqrt{-1}\partial \bar\partial \varphi\Lambda u, u\rangle _h=0$. 
\end{rem}
\end{say}

\section{Proof of the main theorem}\label{sec3} 
In this section, we prove Theorem \ref{II}. 
So, we freely use the notation in Theorem \ref{II}. 
Let $W\Subset Y$ be any Stein open subset. 
We put $V=f^{-1}(W)$. Then $V$ is a holomorphically 
convex weakly $1$-complete K\"ahler manifold (see Remark \ref{remark31} below). 
Let $\mathcal F$ be a coherent sheaf on $V$. 
Then 
$$
f^*: H^q(V, \mathcal F)\to \Gamma (W, R^qf_*\mathcal F)
$$ 
is an isomorphism of topological vector spaces for every $q$ 
(cf.~\cite[Lemma II.1.]{prill}). 
In particular, 
$H^q(V, \mathcal F)$ is a separated topological space. 
Therefore, to prove Theorem \ref{II}, it is sufficient to show that 
$$
\times s: H^q(V, K_V\otimes E\otimes F\otimes \mathcal J(h_F))
\to H^q(V, K_V\otimes E\otimes F\otimes 
\mathcal J(h_F)\otimes L)
$$ is 
injective for every $q\geq 0$. 

\begin{rem}\label{remark31}
A weakly $1$-complete manifold $X$ 
is called a {\em{weakly pseudoconvex 
manifold}} in \cite{dem}. A weakly $1$-complete manifold is a complex 
manifold equipped with a {\em{smooth}} plurisubharmonic 
exhaustion function. More explicitly, there exists a 
smooth plurisubharmonic function $\varphi$ on $X$ such 
that $X_c=\{x\in X | \varphi(x)<c\}$ is relatively compact in 
$X$ for every $c$. 
\end{rem}

We define bounded smooth functions from the 
given nonzero holomorphic section $s$ of $L$. 
\begin{defn}
Take a smooth plurisubharmonic exhaustion function $\varphi$ 
on 
$V$. Without loss of generality, we can assume that 
$\min_{x\in V} \varphi(x)=0$. 
Let $s$ be a holomorphic section of $L$. 
Let $|s|$ be the pointwise norm of $s$ with respect to 
the fiber metric $h_L$. 
Let $\lambda:[0, \infty)\to [0, \infty)$ be a 
smooth convex increasing function such that 
$5|s|^2<e^{\lambda(\varphi)}$. 
Thus, $|s|^2_{\lambda(\varphi)}<\frac{1}{5}<\frac{1}{4}$, where 
$|s|_{\lambda(\varphi)}$ is the pointwise norm of 
$s$ with respect to the fiber metric $h_Le^{-\lambda(\varphi)}$. 
We put $\mu(x):=\lambda(x)+x$. Obviously, $\mu$ is also a smooth 
convex increasing function. 
\end{defn}

\begin{defn}
We put 
$
\chi(t)=t-\log (-t)
$
for $t<0$. 
We define 
\begin{eqnarray*}
\sigma _{\varepsilon, \lambda}&=&\log (|s|^2_{\lambda(\varphi)}
+\varepsilon), \ \text{and}\\
\eta_{\varepsilon, \lambda}&=
&\frac{1}{\epsilon}-\chi(\sigma_{\varepsilon, \lambda})\\
&=& -\log (|s|^2_{\lambda(\varphi)}
+\varepsilon)+\log (-\log (|s|^2_{\lambda(\varphi)}
+\varepsilon))+\frac{1}{\varepsilon} 
\end{eqnarray*} 
for a small positive real number $\varepsilon$. 
We can also define $\sigma_{\varepsilon, \mu}$ and 
$\eta_{\varepsilon, \mu}$ similarly. 
Note that $\eta_{\varepsilon, \lambda}$ and 
$\eta_{\varepsilon, \mu}$ are 
smooth bounded functions on $V$ with 
$\eta_{\varepsilon, \mu}\geq 
\eta_{\varepsilon, \lambda}>\frac{1}{\varepsilon}$. 
The subscripts $\lambda$, $\mu$, and 
$\varepsilon$ might be 
suppressed if there is no danger of confusion. 
\end{defn}

We note the following obvious remark 
before we start various calculations. 

\begin{rem}\label{rem34} 
We note that $e<2\sqrt{2}$. 
Thus, $\frac{3}{2}\log 2>1$. 
Therefore, $\sigma_{\varepsilon, \mu}\leq 
\sigma_{\varepsilon, \lambda}
<2\log \frac{1}{2}
<-\frac{4}{3}$ if $\varepsilon$ is small since 
$|s|^2_{\mu(\varphi)}+\varepsilon
\leq |s|^2_{\lambda(\varphi)}+\varepsilon<\frac{1}{4}$ 
by $|s|^2_{\lambda(\varphi)}<\frac{1}{5}$. 
Of course, $\log(-\sigma_{\varepsilon,\mu})\geq 
\log (-\sigma_{\varepsilon,\lambda})>\log 
\frac{4}{3}>0$. We have 
$\chi'(t)=1-\frac{1}{t}$ and $\chi''(t)=\frac{1}{t^2}$. 
Thus, $1<\chi'
(\sigma_{\varepsilon, \mu})\leq 
\chi'(\sigma_{\varepsilon, 
\lambda})=
1+\frac{1}{(-\sigma_{\varepsilon, \lambda})}<\frac{7}{4}$. 
We also note that 
$-\chi (\sigma_{\varepsilon, \mu})\geq -\chi (\sigma_{\varepsilon, \lambda})
\geq \frac{4}{3}+\log \frac{4}{3}>0$. 
\end{rem}

\begin{say}[Basic calculations]\label{cal} 
We calculate various differentials of 
$\eta_{\varepsilon, \lambda}$. 
The same arguments work for $\eta_{\varepsilon, \mu}$.  
\begin{defn}
Let $u\in C^{p,q}(V,L)$ and $v\in C^{r.s}(V,L)$. 
We define 
$$
\{u, v\}_{\lambda(\varphi)}=u\wedge H_Le^{-\lambda(\varphi)}
\bar v,  
$$ 
where 
$H_L$ is the local matrix representation of $h_L$. 
\end{defn}

We have 
$$\partial \sigma_{\varepsilon,\lambda}=
\frac{
\{D's, s\}_{\lambda(\varphi)}}{|s|^2_{\lambda(\varphi)}
+\varepsilon}. $$
In the above equation, $D'$ is the $(1,0)$ part of the 
Chern connection of $L'=(L, h_Le^{-\lambda(\varphi)})$, 
that is, $D'=D'_{(L, h_Le^{-\lambda(\varphi)})}$. 
Thus, $\Theta (L')=
\Theta (L)+\partial \bar\partial 
\lambda(\varphi)$. 
We obtain the following equation by the direct computation. 
\begin{align*}
\sqrt{-1}\partial\bar\partial
\sigma_{\varepsilon,\lambda}=&
-\frac{
\{\sqrt{-1}\Theta(L')s, s\}_{\lambda(\varphi)}}{|s|^2_{\lambda(\varphi)}
+\varepsilon}\\
&+\frac{\sqrt{-1}
\{D's, D's\}_{\lambda(\varphi)}}{|s|^2_{\lambda(\varphi)}
+\varepsilon}
-\frac{\sqrt{-1}
\{D's, s\}_{\lambda(\varphi)}
\wedge \{s, D's\}_{\lambda(\varphi)}}
{(|s|^2_{\lambda(\varphi)}
+\varepsilon)^2}, 
\end{align*}
where $L'=(L, h_Le^{-\lambda(\varphi)})$. 
Since $L$ is a line bundle, 
we have 
\begin{eqnarray*}
\sqrt{-1}
\{D's, D's\}_{\lambda(\varphi)}|s|^2_{\lambda(\varphi)} 
=\sqrt{-1}
\{D's, s\}_{\lambda(\varphi)}
\wedge \{s, D's\}_{\lambda(\varphi)}. 
\end{eqnarray*}
Substituting it 
into the above equation, we obtain 
\begin{align*}
\sqrt{-1}\partial\bar\partial
\sigma_{\varepsilon,\lambda}= & 
\, \frac{\varepsilon}{|s|^2_{\lambda(\varphi)}
(|s|^2_{\lambda(\varphi)}+\varepsilon)^2}
\sqrt{-1}\{D's, s\}_{\lambda(\varphi)}
\wedge \{s, D's\}_{\lambda(\varphi)}\\ &-\frac{
\{\sqrt{-1}\Theta(L')s, s\}_{\lambda(\varphi)}}{|s|^2_{\lambda(\varphi)}
+\varepsilon}\\ 
=&\, 
\frac{\varepsilon}{|s|^2_{\lambda(\varphi)}}\sqrt{-1}
\partial \sigma_{\varepsilon,\lambda}\wedge\bar\partial 
\sigma_{\varepsilon,\lambda}-\frac{
\{\sqrt{-1}\Theta(L')s, s\}_{\lambda(\varphi)}}{|s|^2_{\lambda(\varphi)}
+\varepsilon}. 
\end{align*}
\begin{lem}
The next equations easily follow from the definition. 
\begin{align*}
\partial\eta_{\varepsilon, \lambda}&=
-\chi'(\sigma_{\varepsilon, \lambda})
\partial\sigma_{\varepsilon, \lambda}, \\ 
\bar\partial\eta_{\varepsilon, \lambda}&=
-\chi'(\sigma_{\varepsilon, \lambda})
\bar\partial\sigma_{\varepsilon, \lambda}, \\ 
\partial \bar\partial \eta_{\varepsilon, \lambda}&=
-\chi''(\sigma _{\varepsilon, \lambda})\partial 
\sigma_{\varepsilon, \lambda}
\wedge\bar\partial\sigma_{\varepsilon, \lambda}
-\chi'(\sigma_{\varepsilon, \lambda})
\partial \bar\partial\sigma_{\varepsilon, \lambda}. 
\end{align*}
\end{lem}
Combining the above equalities, we have 
\begin{align*}
-\sqrt{-1}\partial \bar\partial \eta_{\varepsilon, \lambda} 
=&\ \chi'(\sigma_{\varepsilon, \lambda})\sqrt{-1} 
\partial \bar\partial\sigma_{\varepsilon, \lambda}+ 
\chi''(\sigma _{\varepsilon, \lambda})\sqrt{-1}\partial 
\sigma_{\varepsilon, \lambda}
\wedge\bar\partial\sigma_{\varepsilon, \lambda}\\ 
=&
\ \frac{\varepsilon
\chi'(\sigma_{\varepsilon, \lambda})}
{|s|^2_{\lambda(\varphi)}}\sqrt{-1}
\partial \sigma_{\varepsilon, \lambda}\wedge\bar\partial 
\sigma_{\varepsilon, \lambda}
+\chi''(\sigma_{\varepsilon, \lambda})\sqrt{-1}\partial 
\sigma_{\varepsilon, \lambda}\wedge\bar\partial 
\sigma_{\varepsilon, \lambda}
\\&\ -
\chi'(\sigma_{\varepsilon, \lambda})
\frac{
\{\sqrt{-1}\Theta(L')s, s\}_{\lambda(\varphi)}}{|s|^2_{\lambda(\varphi)}
+\varepsilon}\\ 
=&\  \left( \frac{\varepsilon}{\chi'
(\sigma_{\varepsilon, \lambda})|s|^2
_{\lambda(\varphi)}}+\frac{\chi''
(\sigma_{\varepsilon, \lambda})}
{\chi'(\sigma_{\varepsilon, \lambda})^2}\right)\sqrt{-1}\partial 
\eta_{\varepsilon, \lambda}
\wedge\bar\partial\eta_{\varepsilon, \lambda}\\
&\ -\chi'(\sigma_{\varepsilon, \lambda})
\frac{
\{\sqrt{-1}\Theta(L')s, s\}_{\lambda(\varphi)}}{|s|^2_{\lambda(\varphi)}
+\varepsilon}. 
\end{align*}
\begin{lem}
We have the following inequality. 
$$\frac{\chi''
(\sigma_{\varepsilon, \lambda})}
{\chi'(\sigma_{\varepsilon, \lambda})^2}
\geq \eta^{-2}_{\varepsilon, \lambda}$$
\end{lem}
\begin{proof}
By the definition, it is easy to see that 
$$
\frac{\chi'
(\sigma_{\varepsilon, \lambda})}
{\chi''(\sigma_{\varepsilon, \lambda})}
=\frac{(1-\frac{1}{\sigma_{\varepsilon, \lambda}})^2}{\frac{1}
{\sigma^2_{\varepsilon, \lambda}}}
=(\sigma_{\varepsilon, \lambda}-1)^2. 
$$ 
On the other hand, $\eta_{\varepsilon, \lambda}=
\frac{1}{\varepsilon}
-\sigma_{\varepsilon, \lambda}+\log(-\sigma
_{\varepsilon, \lambda})>1-\sigma_{\varepsilon, \lambda}$. 
Thus, we obtain the desired inequality. 
\end{proof}
\end{say}
By the above calculations, we obtain a very important inequality. 
\begin{prop}\label{34}  
Under the curvature conditions  
$$
\sqrt{-1}(\Theta(E)+\xId_E\otimes \Theta(F))\geq_{\Nak}0
$$ 
on $X\setminus Z$ and 
$$\sqrt{-1}(\Theta(E)+\xId_E\otimes \Theta(F)-
\varepsilon_0 
\xId_E\otimes \Theta(L))\geq_{\Nak}0$$ on $X\setminus Z$ 
for some positive 
real number $\varepsilon_0$, 
we have a small positive real number $\varepsilon_1$ such 
that 
$$
\sqrt{-1}(\eta\Theta_{(E\otimes F,h_Eh_Fe^{-\mu(\varphi)})}-
\xId_E\otimes \partial \bar\partial \eta)\geq_{\Nak} 
\sqrt{-1}(\xId_E\otimes \eta^{-2}\partial \eta\wedge\bar\partial \eta)  
$$ 
holds on $V\setminus Z$ for $0<\varepsilon<\varepsilon_1$, 
where $\eta=\eta_{\varepsilon, \lambda}$ or $\eta_{\varepsilon, 
\mu}$. 
\end{prop}

\begin{proof} 
By the definitions of $\lambda$ and $\mu$, 
$\partial \bar\partial \mu(\varphi)=
\partial \bar\partial \lambda(\varphi)+
\partial \bar\partial \varphi$, and 
$\lambda(\varphi)$ and $\mu(\varphi)$ are 
plurisubharmonic. Note that 
\begin{align*}
\Theta_{(E\otimes F, 
h_Eh_Fe^{-\mu(\varphi)})}=&\ \Theta(E)+\xId_E\otimes 
\Theta(F)+\xId _E\otimes\partial 
\bar\partial \mu(\varphi)\\
=&\ \Theta(E)+\xId_E\otimes \Theta(F)+\xId _E\otimes 
\partial \bar\partial \lambda(\varphi)
+\xId_E\otimes\partial \bar\partial \varphi, \\
\Theta_{(L, h_Le^{-\mu(\varphi)})}=&\ \Theta(L)
+\partial \bar\partial \mu(\varphi), \ \text{and}\\
\Theta_{(L,h_Le^{-\lambda(\varphi)})}=&\ \Theta(L)+
\partial \bar\partial \lambda(\varphi). 
\end{align*}
We also note that 
$$
0\leq \chi'(\sigma _{\varepsilon, \lambda})\frac{|s|^2_{\lambda (\varphi)}}
{|s|^2_{\lambda(\varphi)}+\varepsilon}<\frac{7}{4}
$$ 
and 
$$
0\leq \chi'(\sigma _{\varepsilon, \mu})\frac{|s|^2_{\mu (\varphi)}}
{|s|^2_{\mu(\varphi)}+\varepsilon}<\frac{7}{4}
$$ 
by Remark \ref{rem34}. 
If $\varepsilon_1$ is small, then 
$$
\eta\geq \max \left\{\frac{7}{4\varepsilon _0}, \frac{7}{4}\right\}
$$  
since $\eta>\frac{1}{\varepsilon}>\frac{1}{\varepsilon_1}$. 
Therefore, 
$$
\sqrt{-1}(\eta\Theta_{(E\otimes F,h_Eh_Fe^{-\mu(\varphi)})}-
\xId_E\otimes \partial \bar\partial \eta)\geq_{\Nak} 
\sqrt{-1}(\xId_E\otimes \eta^{-2}\partial \eta\wedge\bar\partial \eta)  
$$  
holds on $V\setminus Z$ for $0<\varepsilon <\varepsilon _1$ where 
$\eta=\eta_{\varepsilon, \lambda}$ or $\eta_{\varepsilon, \mu}$. 
\end{proof}
We note that 
we need no assumptions on the sign of 
the curvature $\sqrt{-1}\Theta (L)$ in Proposition \ref{34}. 
It is a very important remark. 

In the next lemma, we obtain the relationship between 
the Chern connections of $(L, h_L)$ and $(L, h_Le^{-\lambda(\varphi)})$. 
\begin{lem}Let $\gamma :[0, \infty)\to \mathbb R$ be any 
smooth $\mathbb R$-valued function. 
Then we have the following equation by the definition of the 
Chern connection. 
\begin{align*}
D'_{(L, h_Le^{-\gamma(\varphi)})}=
&\ (H_Le^{-\gamma(\varphi)})^{-1}
\partial (H_Le^{-\gamma(\varphi)} \cdot\ )\\ 
=&\ \partial +\partial \log (H_Le^{-\gamma(\varphi)})\wedge\cdot\\
=&\ \partial +\partial \log H_L\wedge\cdot \,-\gamma'(\varphi)\partial 
\varphi\wedge\cdot\\ 
=&\ D'_{(L, h_L)}-\gamma'(\varphi)\partial \varphi\wedge\cdot\ . 
\end{align*}
We note that $H_L=\overline H_L$ since $L$ is a line bundle. 
\end{lem} 

\begin{say}[Complete K\"ahler metrics] 
There exists a complete K\"ahler metric $g$ on $V$ 
since $V$ is weakly $1$-complete. Let $\omega$ be 
the fundamental form of $g$. 
We note the following well-known lemma (cf.~\cite[Lemma 5]{dem2}). 
\begin{lem}
There exists a quasi-psh function $\psi$ on $X$ such that 
$\psi=-\infty$ on $Z$ with logarithmic poles along $Z$ and $\psi$ 
is smooth outside $Z$. 
\end{lem}
Without loss of generality, we can assume that $\psi< -e$ 
on $V\Subset X$. 
We put $\widetilde {\psi}
=\frac{1}{\log (-\psi)}$. Then $\widetilde{\psi}$ 
is a quasi-psh function 
on $V$ and 
$\widetilde \psi< 1$. 
Thus, we can take a positive constant 
$\alpha$ such that $\sqrt{-1} \partial \bar\partial 
\widetilde \psi+\alpha \omega> 0$ 
on $V\setminus Z$. 
Let $g'$ be the K\"ahler metric on $V\setminus Z$ whose 
fundamental form is $\omega'=\omega+(\sqrt{-1} \partial \bar\partial 
\widetilde \psi+\alpha \omega)$. We note that we can check that 
$$
\omega'\geq \sqrt{-1}\partial (\log (\log (-\psi)))\wedge 
\bar\partial (\log (\log (-\psi)))  
$$ if we choose $\alpha \gg 0$. 
It is because 
$$\partial \bar\partial \widetilde \psi=
2\frac{
\frac{-\partial \psi}{-\psi}
\wedge \frac{-\bar\partial \psi}{-\psi}
}{(\log (-\psi))^3}
+
\frac{
\frac{\partial \bar\partial \psi}{-\psi}
}{(\log (-\psi))^2}
+
\frac{
\frac{-\partial \psi \wedge (-\bar\partial \psi)}
{(-\psi)^2}}
{(\log (-\psi))^2}, 
$$ and 
$$
\partial (\log(\log (-\psi)))=\frac{
\frac{-\partial \psi}{-\psi}}{\log (-\psi)}. 
$$ 
Therefore, $g'$ is a {\em{complete}} K\"ahler metric 
on $V\setminus Z$ by 
Hopf--Rinow because 
$\log (\log (-\psi))$ tends to $+\infty$ on $Z$. 
For similar arguments, see \cite[Section 3]{f2}. 
We fix these K\"ahler metrics throughout this proof. 
\end{say}

\begin{say}[Key Results] 
The following three propositions are the heart of 
the proof of Theorem \ref{II}.  

\begin{prop}\label{314}
For every $u\in \mathcal H^{n,q}(V\setminus Z, (E\otimes F, h_Eh_F
e^{-\mu(\varphi)}))_{g'}$, 
$(\bar\partial \eta)^*u=0$ 
for $\eta=\eta_{\varepsilon, \lambda}$ and 
$\eta_{\varepsilon, \mu}$. 
This implies that $\partial \eta\wedge 
\ast u=0$ for $\eta=\eta_{\varepsilon, \lambda}$ 
and $\eta_{\varepsilon, \mu}$. Thus, we obtain 
$D'_{(L, h_Le^{-\lambda(\varphi)})}s\wedge \ast u=0$ 
and $D'_{(L, h_Le^{-\mu(\varphi)})}s\wedge \ast u=0$. 
\end{prop}

\begin{proof}
The definition of $\mathcal H^{n,q}(V\setminus Z, (E\otimes F, h_E
h_Fe^{-\mu(\varphi)}))_{g'}$ 
implies that $\bar\partial u=0$ and 
$D''^*_{(E\otimes F,h_Eh_Fe^{-\mu(\varphi)})}u=0$. 
By Propositions \ref{twist}, \ref{24}, and 
\ref{34}, we have 
$$
0\geq -\|\sqrt{\eta}{D'}^*u\|^2\geq \lla \sqrt{-1}
\eta^{-2}\partial \eta\wedge\bar\partial\eta \Lambda u, 
u\rra\geq 0. 
$$ 
Thus, we have $(\bar\partial \eta)^*u=0$ (cf.~Lemma \ref{26}). 
Therefore, we obtain 
$\partial \eta\wedge \ast u=0$ because 
$(\bar\partial \eta)^*=\ast \partial \eta\wedge \ast$ 
by Lemma \ref{lem-adj}. 
By the definition of $\eta$, we obtain 
${D'}_{(L, h_Le^{-\lambda(\varphi)})}s\wedge\ast u=0$ 
(resp.~${D'}_{(L, h_Le^{-\mu(\varphi)})}s\wedge\ast u=0$) 
if $\eta=\eta_{\varepsilon,\lambda}$ (resp.~$\eta=\eta
_{\varepsilon, \mu}$).  
\end{proof}

\begin{prop}\label{315}
If $D'_{(L, h_Le^{-\lambda(\varphi)})}s\wedge\ast u=0$ and 
$D'_{(L, h_Le^{-\mu(\varphi)})}s\wedge\ast u=0$, 
then $D'_{(L, h_L)}s\wedge \ast u=0$ and $\partial \varphi 
\wedge \ast u=0$. 
Therefore, 
$D'_{(L, h_Le^{-\nu(\varphi)})}s\wedge \ast u=0$ for 
any smooth $\mathbb R$-valued function $\nu$ 
defined on $[0,\infty)$. 
\end{prop}

\begin{proof}
We note that $D'_{(L, h_Le^{-\lambda(\varphi)})}=
D'_{(L, h_L)}-\lambda'(\varphi)\partial \varphi\wedge\cdot\ $ 
and 
$$D'_{(L, h_Le^{-\mu(\varphi)})}=
D'_{(L, h_L)}-\lambda'(\varphi)\partial \varphi\wedge\cdot\,-
\partial \varphi\wedge\cdot\ $$ since 
$\mu(x)=\lambda(x)+x$.  
\end{proof}

\begin{prop}\label{211} 
Assume that $D'_{(L, h_Le^{-\nu(\varphi)})}s\wedge *u=0$ for a 
smooth $\mathbb R$-valued function $\nu$ defined on $[0, \infty)$.   
If $D''^*_{(E\otimes F, h_Eh_Fe^
{-\mu(\varphi)})}u=0$, then we obtain 
$$D''^*_{(E\otimes F\otimes L, h_Eh_Fh_Le^
{-\mu(\varphi)-\nu(\varphi)})}(su)=0. $$ 
\end{prop}
\begin{proof} Let $H_E$ (resp.~$H_F$) be the local matrix representation of 
$h_E$ (resp.~$h_F$). 
The condition $D''^*_{(E\otimes F, h_Eh_Fe^{-\mu(\varphi)})}u=0$ 
implies 
that 
$$\bar\partial (e^{-\mu(\varphi)}H_EH_F\overline {\ast u})=0.$$ 
To prove $D''^*_{(E\otimes F\otimes L, h_Eh_Fh_Le^{-\mu(\varphi)
-\nu(\varphi)})}(su)=0$, 
it is sufficient to check that 
$\bar\partial (H_EH_Fe^{-\mu(\varphi)
-\nu(\varphi)}H_L\overline {\ast su})=0$. 
We note that 
$$
\bar\partial (H_EH_Fe^{-\mu(\varphi)-\nu
(\varphi)}H_L\overline {\ast su})=
\bar\partial (H_L\overline {s}e^{-\nu(\varphi)})\wedge e^
{-\mu(\varphi)}H_EH_F\overline 
{\ast u}
$$ by the above condition. 
The right hand side is zero since 
$D'_{(L, h_Le^{-\nu(\varphi)})}s\wedge \ast u=0$. 
\end{proof}
\end{say}

The next theorem is a key result. 

\begin{thm}[{cf.~\cite[Proposition 3.1]{ohsawa}}]\label{310}
For any smooth $\mathbb R$-valued function $\nu$ 
defined on $[0,\infty)$ such that 
$\nu\geq C$ for some constant $C$, 
$$s\mathcal H^{n,q}(V\setminus Z, (E\otimes F, h_Eh_F
e^{-\mu(\varphi)}))_{g'}$$ is contained 
in $$\mathcal H^{n,q}(V\setminus Z, (E\otimes F\otimes L, h_Eh_Fh_L
e^{-\mu
(\varphi)-\lambda(\varphi)-\nu(\varphi)}))_{g'}$$ for every $q$. 
\end{thm}
\begin{proof} 
Let $u\in \mathcal H^{n,q}(V\setminus Z, (E\otimes F, h_Eh_F
e^{-\mu(\varphi)}))_{g'}$. 
Then it is obvious that 
$su\in L^{n,q}_{(2)}(E\otimes F\otimes L, h_Eh_Fh_Le
^{-\mu(\varphi)-\lambda(\varphi)-\nu(\varphi)})$ 
because $|s|^2_{\lambda(\varphi)}<\frac{1}{5}$ 
and $0<e^{-\nu(\varphi)}\leq e^{-C}$. 
So, the claim is a direct consequence of Propositions \ref{314}, 
\ref{315}, and \ref{211}. 
Note that $\bar\partial (su)=0$ for 
$u\in \mathcal H^{n,q}(V\setminus Z, (E\otimes F, h_Eh_F
e^{-\mu(\varphi)}))_{g'}$ 
since $s$ is holomorphic and 
$\bar\partial u=0$. 
\end{proof}

\begin{say}[Cohomology groups]
Before we start the proof of the main theorem:~Theorem \ref{II}, 
we represent the cohomology groups on $V$ 
by the objects on $V\setminus Z$. 
\begin{defn}[Space of locally square integrable forms] 
We define the space of locally (in $V$) square 
integrable $E\otimes F$-valued $(n,q)$-forms on $V\setminus Z$. 
It is denoted by $L^{n,q}_{loc, V}(V
\setminus Z, E\otimes F)$ or 
$L^{n,q}_{loc, V}(V
\setminus Z, (E\otimes F, h_Eh_F))$. 
The vector space $L^{n,q}_{loc, V}(V
\setminus Z, E\otimes F)$ is spanned by $(n,q)$-forms 
$u$ on $V\setminus Z$ 
with measurable coefficients such that 
$$
\int_{U}|u|^{2}_{g',h_Eh_F}dV_{\omega'}<\infty 
$$ 
for every $U\Subset V$ (not $U\Subset V\setminus Z$), 
where $|\cdot|_{g', h_Eh_F}$ is the pointwise norm 
with respect to $g'$ and $h_Eh_F$.
We note the following obvious 
remark. Let $h:V\to (0, \infty)$ be a smooth 
positive function. Then 
$L^{n,q}_{loc, V}(V
\setminus Z, (E\otimes F, h_Eh_F))=
L^{n,q}_{loc, V}(V
\setminus Z, (E\otimes F, h_Eh_Fh))$. 
We can define $L^{n,q}_{loc, V}(V
\setminus Z, E\otimes F\otimes L)$ similarly. 
\end{defn}

The next lemma is essentially the same as \cite[Claim 1]{f2}, 
which is more or less known to 
experts (cf.~\cite[Proposition 4.6]{take}). 

\begin{lem}\label{320}The following isomorphism holds. 
\begin{eqnarray*}
&&H^q(V, K_V\otimes E\otimes F\otimes \mathcal J(h_F))\\
&\simeq& H^{n,q}_{loc, V}(V\setminus Z, E\otimes F)_{g'}\\
&:=&\frac{\Ker \bar\partial \cap 
L^{n,q}_{loc, V}(V\setminus Z,E\otimes F)}
{L^{n,q}_{loc, V}(V\setminus Z,E\otimes F)
\cap \bar\partial L^{n,q-1}_{loc, V}(V\setminus Z,E\otimes F)}. 
\end{eqnarray*}
\end{lem}
\begin{proof}[Sketch of the proof] 
Let $V=\bigcup _{i\in I}U_i$ be a locally finite Stein cover of $V$ 
such that each $U_i$ is sufficiently small and $U_i\Subset V$. 
We denote this cover by $\mathcal U=\{U_i\}_{i\in I}$. 
By Cartan and Leray, we obtain 
$$H^q(V, K_V\otimes E\otimes F\otimes \mathcal J(h_F))\simeq 
{\check{H}}^q(\mathcal U, K_V\otimes E\otimes F\otimes \mathcal J(h_F)),$$ 
where the right hand side is the $\check{\rm C}$ech 
cohomology group calculated 
by $\mathcal U$. 
By using a partition of unity $\{\rho_i\}_{i \in I}$ 
associated to $\mathcal U$, we can construct a homomorphism  
$$
\rho: {\check{H}}^q({\mathcal U}, K_X\otimes E\otimes F\otimes \mathcal J(h_F))
\to 
H^{n,q}_{loc, V}(V\setminus Z, E\otimes F)_{g'}.   
$$ 
See Remark \ref{rem215} below. 
We can check that $\rho$ is an isomorphism. 
Note the following 
facts:~(a) The open set $U_{i_0}\cap\cdots \cap U_{i_k}$ 
is Stein. So, $U_{i_0}\cap\cdots \cap U_{i_k}\setminus Z$ is a complete 
K\"ahler manifold (cf.~\cite[Th\'eor\`eme 0.2]{dem-est}). 
Therefore, $E\otimes F$-valued $\bar\partial$-equations can be 
solved with suitable 
$L^2$ estimates on $U_{i_0}\cap\cdots \cap U_{i_k}\setminus Z$ 
by Lemma \ref{3.2} 
below. 
(b) Let $U$ be an open subset of $V$. 
An $E\otimes F$-valued holomorphic $(n,0)$-form on $U\setminus Z$ 
with 
a finite $L^2$ norm can be extended to an $E\otimes F$-valued holomorphic 
$(n,0)$-form on $U$ (cf.~Remark \ref{rem215}). 
\end{proof}

\begin{rem}[{cf.~\cite[Lemme 3.3]{dem-est}}]\label{rem215}
Let $u$ be an $E\otimes F$-valued $(n,q)$-form on $V\setminus 
Z$ with measurable coefficients. 
Then, we have 
$|u|^2_{g', h_Eh_F}dV_{\omega'}\leq |u|^2_{g, h_Eh_F}dV_{\omega}$, where 
$|u|_{g', h_Eh_F}$ (resp.~$|u|_{g, h_Eh_F}$) is the pointwise norm induced 
by $g'$ and $h_Eh_F$ (resp.~$g$ and $h_Eh_F$) since 
$g'>g$ on $V\setminus Z$. 
If $u$ is an $E\otimes F$-valued $(n,0)$-form, then 
$|u|^2_{g', h_Eh_F}dV_{\omega'}= |u|^2_{g, h_Eh_F}dV_{\omega}$. 
\end{rem}

The following lemma is \cite[Lemma 3.2]{f2}, which is a reformulation of the 
classical $L^2$-estimates for $\bar\partial$-equations for our purpose. 

\begin{lem}[$L^2$-estimates for $\bar\partial$-equations on 
complete K\"ahler manifolds]\label{3.2}
Let $U$ be a sufficiently small Stein open set of $V$. If 
$u\in L^{n,q}_{(2)}(U\setminus Z, E\otimes F)_{g', h_Eh_F}$ 
with $\bar\partial u=0$, then 
there exists 
$v\in L^{n,q-1}_{(2)}(U\setminus Z, E\otimes F)_{g', h_Eh_F}$ such that 
$\bar\partial v =u$. 
Moreover, there exists a positive constant $C$ independent 
of $u$ such that 
$$
\int _{U\setminus Z} |v|^2_{g', h_Eh_F}\leq C\int _{U\setminus Z}
|u|^2_{g', h_Eh_F}. 
$$
We note that $g'$ is not a complete K\"ahler metric on $U\setminus Z$ but 
$U\setminus Z$ is a complete K\"ahler 
manifold {\em{(}}cf.~\cite[Th\'eor\`eme 0.2]{dem-est}{\em{)}}. 
\end{lem}

By the same arguments, the isomorphism in Lemma \ref{320} 
holds even when we replace $(E\otimes F, h_Eh_F)$ with 
$(E\otimes F\otimes L, h_Eh_Fh_L)$. 
\end{say}

\begin{say}[{Proof of the main theorem:~Theorem \ref{II}}] 
Let us start the proof of Theorem \ref{II} (cf.~\cite{ohsawa}). 
\begin{proof}[{Proof of {\em{Theorem \ref{II}}}}]
Let $u$ be any $\bar\partial$-closed locally 
square integrable $E\otimes F$-valued $(n,q)$-form 
on $V\setminus Z$ such that $su=\bar\partial v$ for 
some $v\in L^{n,q-1}_{loc, V}
(V\setminus Z, E\otimes F\otimes L)$. 
We choose $\lambda$ such that $|s|^2_{\lambda(\varphi)}<\frac{1}{5}$ and 
$u\in L^{n,q}_{(2)}(V\setminus Z, (E\otimes F, 
h_Eh_Fe^{-\lambda(\varphi)}))$. 
Since $\mu(\varphi)=\lambda(\varphi)+\varphi$, 
$u\in L^{n,q}_{(2)}(V\setminus  Z, (E\otimes F, 
h_Eh_Fe^{-\mu(\varphi)}))$. 
By choosing $\nu$ suitably, we can assume that 
$v\in L^{n,q-1}_{(2)}(V\setminus Z, (E\otimes F\otimes L, 
h_Eh_Fh_Le^{-\mu(\varphi)-\nu(\varphi)}))$. 
In particular, $v\in L^{n,q-1}_{(2)}(V\setminus Z, (E\otimes 
F\otimes L, h_Eh_Fh_L
e^{-\mu(\varphi)-\lambda(\varphi)-\nu(\varphi)}))$. 
Let $Pu$ be the orthogonal projection of $u$ 
to $\mathcal H^{n,q}(V\setminus Z, (E\otimes F, 
h_Eh_Fe^{-\mu(\varphi)}))_{g'}$. 
We note that 
\begin{eqnarray*}
&&L^{n,q}_{(2)}(V\setminus Z, 
(E\otimes F, h_Eh_Fe^{-\mu(\varphi)}))\\
&=&\overline{\xIm \bar\partial}
\oplus \mathcal H^{n,q}(V\setminus Z, (E\otimes F, h_Eh_F
e^{-\mu(\varphi)}))_{g'}\oplus 
\overline{\xIm D''^*_{(E\otimes F, 
h_Eh_Fe^{-\mu(\varphi)})}},
\end{eqnarray*}
and 
\begin{eqnarray*}
\Ker \bar\partial&=&
\overline{\xIm \bar\partial}
\oplus \mathcal H^{n,q}(V\setminus Z, (E\otimes F, h_Eh_F
e^{-\mu(\varphi)}))_{g'}.  
\end{eqnarray*}
Here, $\overline{\xIm \bar\partial}$ 
(resp.~$\overline{\xIm D''^*_{(E\otimes F, 
h_Eh_Fe^{-\mu(\varphi)})}}$) denotes 
the closure of $\bar\partial C_0^{n, q-1}(V\setminus Z, 
E\otimes F)$ 
(resp.~$D''^*_{(E\otimes F, 
h_Eh_Fe^{-\mu(\varphi)})}C_0^{n,q+1}(V\setminus Z, E\otimes 
F)$) in 
$L^{n. q}_{(2)}(V\setminus Z, (E\otimes F, 
h_Eh_Fe^{-\mu(\varphi)}))$. 
Note that the fixed K\"ahler metric $g'$ is complete. 
Therefore, $u-Pu$ is in the closure of the image of $\bar\partial$. 
Thus, so is $s(u-Pu)$ since $s$ is holomorphic. 
On the other hand, 
$$sPu\in \mathcal H^{n,q}(V\setminus Z, 
(E\otimes F\otimes L, 
h_Eh_Fh_Le^{-\mu(\varphi)-\lambda(\varphi)
-\nu(\varphi)}))_{g'}$$ by Theorem \ref{310}. 
So, $sPu$ coincides with the orthogonal projection of $su$ to 
$\mathcal H^{n,q}(V\setminus Z, (E\otimes F\otimes L, 
h_Eh_Fh_L
e^{-\mu(\varphi)-\lambda(\varphi)-\nu(\varphi)}))_{g'}$, which 
must be equal to zero since $v\in 
L^{n,q-1}_{(2)}(V\setminus Z, (E\otimes F\otimes L, 
h_Eh_Fh_Le^{-\mu(\varphi)-\lambda(\varphi)-\nu(\varphi)}))$. 
Therefore, $Pu=0$. 
Since $H^q(V, K_V\otimes 
E\otimes F\otimes \mathcal J(h_F))$ 
is a separated topological 
vector space (cf.~\cite[Lemma II.1.]{prill}) 
and $u\in \overline{\xIm \bar\partial}$, 
there exists $w\in L^{n,q-1}_{loc, V}(V\setminus Z,E\otimes 
F)$ such that $u=\bar\partial w$ (cf.~\cite[Proposition 4.6]{take} 
and \cite[Claim 1]{f2}). 
This means that 
$u$ represents zero in $H^q(V, K_V\otimes 
E\otimes F\otimes \mathcal J(h_F))$.  
\end{proof}
\end{say}

\section{Corollaries and applications}\label{sec-p}

In this section, we discuss the proofs of 
corollaries in Section \ref{int} and 
some applications of Theorem \ref{II}. 

First, we give a proof of Corollary \ref{to}, which 
is obvious if we apply Theorem \ref{II} 
for $L=\mathcal O_X\simeq f^*\mathcal O_Y$.  

\begin{proof}[{Proof of {\em{Corollary \ref{to}}}}]
The statement is local. 
So, we can assume that $Y$ is Stein. 
Let $s\in H^0(Y,\mathcal O_Y)$ be 
an arbitrary nonzero section. 
By Theorem \ref{II}, 
$$
\times s: R^qf_*(K_X\otimes E\otimes 
F\otimes \mathcal J(h_F))\to 
R^qf_*(K_X\otimes E\otimes 
F\otimes \mathcal J(h_F))
$$ is injective for every $q\geq 0$. 
Thus, $R^qf_*(K_X\otimes E\otimes F\otimes \mathcal J(h_F))$ 
is torsion-free for every $q\geq 0$. 
\end{proof}

The following proposition is a slight generalization of 
Theorem \ref{II}. 
\begin{prop}\label{326-ne}
In {\em{Theorem \ref{II}}}, we can weaken the 
assumption that $X$ is K\"ahler as follows. 
For any point $P\in Y$, there exist 
an open neighborhood $U$ of $P$ and a 
proper bimeromorphic morphism 
$g:W\to V:=f^{-1}(U)$ from a K\"ahler manifold $W$. 
\end{prop}
\begin{proof}[Sketch of the proof]
We put $Z'=g^{-1}(Z\cap V)\subset W$. 
We can apply Corollary \ref{to} to $g:W\to V$, $Z'$, 
$(g^*E, g^*h_E)$, and $(g^*F, g^*h_F)$. 
Then we obtain $R^{q}g_*(K_W\otimes g^*E\otimes 
g^*F\otimes \mathcal J(g^*h_F))=0$ for every $q>0$ and 
it is well known (and easy to check) that 
$$g_*(K_W\otimes g^*E\otimes g^*F\otimes \mathcal J(g^*h_F))
\simeq K_V\otimes E\otimes F\otimes \mathcal J(h_F)$$ 
(cf.~\cite[(5.8) Proposition]{dem}). 
Therefore, by Leray's spectral sequence, 
$$R^q(f\circ g)_*(K_W\otimes g^*E\otimes g^*F\otimes \mathcal J(g^*h_F))
\simeq R^qf_*(K_V\otimes E\otimes F\otimes \mathcal J(h_F))$$ 
for every $q\geq 0$. 
Apply Theorem \ref{II} to $f\circ g:W\to U$, $Z'$, $(g^*E, g^*h_E)$, 
$(g^*L, g^*h_L)$, and $(g^*F, g^*h_F)$. 
Then we obtain that 
$$\times s:R^qf_*(K_X\otimes E\otimes F\otimes \mathcal J(h_F))
\to R^qf_*(K_X\otimes E\otimes F\otimes \mathcal J(h_F)\otimes 
L)$$ 
is injective for every $q\geq 0$ by 
the above isomorphisms.  
\end{proof}

The next result is related to the main 
theorem in \cite{levine}. 

\begin{cor}\label{loc}
Let $f:X\to \Delta$ be a smooth 
{\em{projective}} surjective 
morphism from 
a K\"ahler manifold $X$ to a 
disk $\Delta$ and let $E$ be a Nakano semi-positive 
vector bundle on $X$. Assume that there 
exists $D\in | K_{X_0}^{\otimes l}|$ such that $\mathcal J_{X_0}(cD)\simeq 
\mathcal O_{X_0}$ for every $0\leq c<1$, 
where $0\in Y$ and $X_0=f^{-1}(0)$. 
Then there exists an open set $U\subset \Delta$ such that 
$0\in U$ and 
$R^qf_*(K_X^{\otimes m}\otimes E)$ 
is locally free on $U$ for every $q\geq 0$ and 
$1\leq m\leq l$. Equivalently, $\dim_{\mathbb C} 
H^q(X_t, K_{X_t}^{\otimes m}
\otimes E)$ is constant for every $q\geq 0$ by the base change theorem, where 
$t\in U$ and $X_t=f^{-1}(t)$. 
\end{cor}
\begin{proof}
In this proof, we shrink $\Delta$ without mentioning it for simplicity of notation. 
Let $s_0\in H^0(X_0, K_{X_0}^{\otimes l})$ 
such that $D=(s_0=0)$. 
By Siu's extension theorem (see \cite[Theorem 0.1]{siu}), 
there exists $s\in H^0(X, K_{X}^{\otimes l})$ such that 
$s|_{X_0}=s_0$. We consider the singular hermitian metric 
$h_{K_X^{\otimes (m-1)}}=(\frac{1}{|s|^2})^{\frac{m-1}{l}}$ of 
$K_X^{\otimes (m-1)}$. 
By the Ohsawa--Takegoshi $L^2$-extension theorem and the 
assumption $\mathcal J_{X_0}(cD)\simeq \mathcal 
O_{X_0}$ for $0\leq c<1$, 
we obtain that $\mathcal J_X(h_{K_X^{\otimes m-1}})
\simeq \mathcal O_X$ in a neighborhood of $X_0$. 
Therefore, we obtain that 
$R^qf_*(K_X^{\otimes m}\otimes E)
=R^qf_*(K_X\otimes K_X^{\otimes (m-1)}\otimes E\otimes 
\mathcal J_X(h_{K_X^{\otimes m-1}}))$ is locally free by 
Corollary \ref{to}.  
\end{proof}

Let us start the proof of 
the Kawamata--Viehweg--Nadel type vanishing theorem:~Corollary 
\ref{kvn-vani}. 

\begin{proof}[Proof of {\em{Corollary \ref{kvn-vani}}}]
Let $P$ be a point of $Y$. 
The problem is local. 
So we repeatedly shrink $Y$ around $P$ without mentioning 
it explicitly. 
Since $\mathcal M$ is $f$-big, 
we have a bimeromorphic map $\Phi:X\dashrightarrow X'\subset Y\times 
\mathbb P^N$ over $Y$. 
By applying Hironaka's Chow lemma (cf.~\cite[Corollary 2]{hironaka}), 
we can construct a bimeromorphic 
map $\varphi:Z\to X$ from a complex manifold such that 
$f\circ \varphi:Z\to Y$ is projective. 
It is easy to see that 
$$
\varphi_*\left(K_Z\otimes \varphi^*E\otimes 
\varphi^*\mathcal L\otimes \mathcal J\left(\frac{1}
{m}\varphi^*D\right)\right)\\
\simeq K_X\otimes E\otimes \mathcal L\otimes 
\mathcal J\left(\frac{1}{m}D\right). 
$$ 
by the definition of $\mathcal J$ (cf.~\cite[Theorem 9.2.33]{lazarsfeld}). 
\begin{claim}
We have 
$$
R^q\varphi_*\left(K_Z\otimes \varphi^*E\otimes 
\varphi^*\mathcal L\otimes 
\mathcal J\left(\frac{1}{m}\varphi^*D\right)\right)=0
$$
for every $q>0$ 
\end{claim}
\begin{proof}[Proof of {\em{Claim}}]
The problem is local. So we can shrink $X$ and 
assume that 
$\mathcal M$ is trivial. In particular, 
$(\varphi^*\mathcal L)^{\otimes m}\simeq \mathcal O_Z(\varphi^*D)$. 
Thus, by Corollary \ref{to}, 
we obtain that 
$R^q\varphi_*(K_Z\otimes \varphi^*E\otimes 
\varphi^*\mathcal L\otimes \mathcal J(\frac{1}{m}\varphi^*D))$ 
is torsion-free for every $q$. 
Thus, $R^q\varphi_*(K_Z\otimes \varphi^*E\otimes 
\varphi^*\mathcal L\otimes \mathcal J(\frac{1}{m}\varphi^*D))=0$ for 
every $q>0$ since $\varphi$ is bimeromorphic. 
\end{proof}
Therefore, by replacing $X$ with $Z$, we can assume that 
$f:X\to Y$ is projective. 
By Kodaira's lemma, 
we can write $\mathcal L^{\otimes a}\simeq \mathcal A\otimes 
\mathcal O_X(G)$ where 
$a$ is a positive integer, $\mathcal A$ is an $f$-ample 
line bundle on $X$, and $G$ is an effective Cartier divisor on $X$. 
Then $\mathcal L^{\otimes (mb+a)}\simeq 
(\mathcal M^{\otimes b}\otimes \mathcal A)\otimes 
\mathcal O_X(G+bD)$. 
Note that $\mathcal M^{\otimes b}\otimes \mathcal A$ is $f$-ample and 
that $\mathcal J(\frac{G+bD}{mb+a})=\mathcal J(\frac{1}{m}D)$ if 
$b\gg 0$. 
Therefore, we can further assume that 
$\mathcal M$ is $f$-ample. 
Then, by Theorem \ref{II}, we can construct inclusions 
$$
R^qf_*(K_X\otimes E\otimes \mathcal L\otimes \mathcal J)
\\ \subset 
R^qf_*(K_X\otimes E\otimes \mathcal L\otimes 
\mathcal J\otimes \mathcal M^{\otimes k})
$$ 
for every $k>0$ with $|\mathcal M^{\otimes k}|\ne \emptyset$. 
Thus, by Serre's vanishing 
theorem (cf.~\cite[p.~25 Remark (4)]{nakayama2}), we obtain 
$R^qf_*(K_X\otimes E\otimes \mathcal L\otimes \mathcal J)=0$
for every $q>0$. 
\end{proof}

\begin{proof}[Proof of {\em{Corollary \ref{kv-vani}}}]
Note that $m\{H\}$ is Cartier. 
We have $m\ulcorner H\urcorner\sim mH+m\{-H\}$. 
We put $E=\mathcal O_X$, $\mathcal L=\mathcal O_X(\ulcorner H\urcorner)$, 
$\mathcal M=\mathcal O_X(mH)$, and 
$D=m\{-H\}$, and apply Corollary \ref{kvn-vani}. 
Then we obtain 
$R^qf_*(K_X\otimes \mathcal O_X(\ulcorner H\urcorner))=0$ for every 
$q>0$. We note that 
$\mathcal J(\{-H\})=\mathcal O_X$ since 
$\Supp \{H\}$ is a normal crossing divisor 
(cf.~\cite[Lemma 9.3.44]{lazarsfeld}). 
\end{proof}

The proof of Corollary \ref{ko-vani} 
is a routine work. 
So we only sketch the proof here. 
For details, 
see, for example, the 
proofs of \cite[Theorem 1.1 (ii)]{fujino-on} and 
\cite[Theorem 6.3 (ii)]{fujino-fundamental}. 

\begin{proof}[Sketch of the proof of {\em{Corollary \ref{ko-vani}}}]
We can repeatedly shrink $Z$ without mentioning it. 
By Hironaka's Chow lemma and Hironaka's flattening theorem 
(cf.~\cite[Corollary 2, Flattening Theorem]{hironaka}), we can assume that 
$g:Y\to Z$ is projective (cf.~the proof of Corollary 
\ref{kvn-vani}). 
By Kodaira's lemma, we can assume that 
$\mathcal N$ is $g$-ample. 
Let $A$ be a general smooth 
sufficiently $g$-ample Cartier divisor on $Y$. 
We put $B=f^*A$. 
We consider the following short exact sequence
\begin{align*}
0&\to K_X\otimes E\otimes \mathcal L\otimes \mathcal J\to 
K_X\otimes E\otimes \mathcal O_X(B)\otimes \mathcal L\otimes \mathcal J
\\ &\to K_B\otimes E|_B\otimes \mathcal L|_B\otimes \mathcal J|_B\to 0. 
\end{align*}
Since $A$ is general, $\mathcal J|_B=\mathcal J(\frac{1}{m}D|_B)$, 
and 
\begin{align*}
0&\to R^qf_*(K_X\otimes E\otimes \mathcal L\otimes \mathcal J)\to 
R^qf_*(K_X\otimes E\otimes \mathcal O_X(B)\otimes \mathcal L\otimes \mathcal J)
\\ &\to R^qf_*(K_B\otimes E|_B\otimes \mathcal L|_B\otimes \mathcal J|_B)\to 0
\end{align*}
is exact for every $q$. 
By taking the long exact sequence, we obtain 
$$
R^pg_*R^qf_*(K_X\otimes E\otimes \mathcal L\otimes \mathcal J)=0
$$
for every $p\geq 2$ and every $q\geq 0$ because 
$A$ is sufficiently $g$-ample and the induction on dimension. 
Then we obtain the following commutative diagram.  
$$
\xymatrix{R^1g_*R^qf_*(K_X\otimes E\otimes \mathcal L\otimes \mathcal J) 
\ar[d]\ar[r]^{\alpha}& 
R^{1+q}(g\circ f)_*(K_X\otimes E\otimes \mathcal L\otimes \mathcal J) 
\ar[d]^{\beta}\\ 
R^1g_*R^qf_*(K_X\otimes E\otimes \mathcal O_X(B)\otimes 
\mathcal L\otimes \mathcal J) \ar[r] &
R^{1+q}(g\circ f)_*(K_X\otimes E\otimes
\mathcal O_X(B)\otimes \mathcal L\otimes \mathcal J)
}
$$
Note that $\alpha$ is injective by the above vanishing result and 
that $\beta$ is injective by Theorem \ref{II}. 
Since $A$ is sufficiently $g$-ample, 
we have $R^1g_*R^qf_*(K_X\otimes E\otimes \mathcal O_X(B)\otimes 
\mathcal L\otimes \mathcal J)=0$. 
Thus, $R^1g_*R^qf_*(K_X\otimes E\otimes \mathcal L\otimes \mathcal J)=0$ for 
every $q\geq 0$. So we finish the proof.  
\end{proof}

\section{Examples:~nef, semi-positive, and semi-ample line 
bundles}\label{app} 

In this section, we collect some examples 
of nef, semi-positive, and semi-ample 
line bundles. 
These examples help us understand our results 
in \cite{f2} and this paper. 
We think that it is very important to understand 
the differences in the notion of semi-ample, semi-positive, and 
nef line bundles. 

First, we recall the following well-known example. 
It implies that there exists a nef line bundle that has no smooth hermitian 
metrics with semi-positive curvature. 

\begin{ex}[{cf.~\cite[Example 1.7]{dps}}]\label{ex1}
Let $C$ be an elliptic curve and let $\mathcal E$ be 
the rank two vector bundle 
on $C$ which is defined by the unique non-splitting 
extension 
$$
0\longrightarrow \mathcal O_C\longrightarrow 
\mathcal E\longrightarrow \mathcal O_C\longrightarrow 0. 
$$
We consider the ruled surface $X:=\mathbb P_C(\mathcal E)$. 
On that surface there is a unique section 
$D:=\mathbb P_C(\mathcal O_C)\subset 
X$ of $X\to C$ such that 
$\mathcal O_D(D)\simeq\mathcal O_D$ 
and $\mathcal O_X(D)\simeq \mathcal O
_{\mathbb P_C(\mathcal E)}(1)$ is a nef line 
bundle (cf.~\cite[V Proposition 2.8]{hartshorne}). 
It is not difficult to see that 
$H^1(X, K_X\otimes \mathcal O_X(2D))
\to H^1(X, K_X\otimes \mathcal O_X(3D))$ is a zero map, 
$H^1(X, K_X\otimes \mathcal O_X(2D))\simeq \mathbb C$, and 
$H^1(X, K_X\otimes \mathcal O_X(3D))\simeq \mathbb C$. 
We note that $K_X\sim \mathcal O_X(-2D)$. 
Therefore, $\mathcal O_X(D)$ has no smooth hermitian metrics with 
semi-positive curvature by Enoki's injectivity theorem 
(see \cite[Theorem 0.2]{enoki}, \cite[Corollary 1.4]{f2}). Note that 
$\kappa (X, \mathcal O_{\mathbb P_C(\mathcal E)}(1))=0$ 
and $\nu(X, \mathcal O_{\mathbb P_C(\mathcal E)}(1))=1$. 
We also note that Koll\'ar's injectivity theorem implies nothing 
since $\mathcal O_{\mathbb P_C(\mathcal E)}(1)$ is not 
semi-ample.
\end{ex}

The next one is an example 
of nef and big line bundles that have no 
smooth hermitian metrics with 
semi-positive curvature.  
The author learned the following construction from Dano Kim. 

\begin{ex}\label{ex52}
We use the same notation as in Example \ref{ex1}. 
Let $P\in C$ be a closed point. 
We put $\mathcal F:=\mathcal E\oplus \mathcal O_C(P)$ 
and $Y:=\mathbb P_C(\mathcal F)$. 
Then it is easy to see that $\mathcal L:=\mathcal O
_{\mathbb P_C(\mathcal F)}(1)$ is nef and big 
(cf.~\cite[Example 6.1.23]{lazarsfeld}). 
Since $\mathcal O_{\mathbb P_C(\mathcal E)}(1)$ 
has no smooth hermitian metrics with 
semi-positive curvature, neither has $\mathcal L$. 
In this case, $H^i(Y, K_Y\otimes \mathcal L^{\otimes k})=0$ for 
$i>0$ and 
every $k\geq 1$ by the Kawamata--Viehweg vanishing theorem. 
\end{ex} 

Let us recall some examples 
of semi-positive line bundles that 
are not semi-ample. 

\begin{ex}
[{cf.~\cite[p.145]{del}}]\label{ex3}
Let $C$ be a smooth projective curve with the genus $g(C)\geq 1$. 
Let $L\in \Pic ^0(C)$ be non-torsion. 
We put $\mathcal E:=\mathcal O_C\oplus 
L$ and $X:=\mathbb P_C(\mathcal E)$. 
Then $\mathcal L:=\mathcal O_{\mathbb P_C(\mathcal 
E)}(1)$ is semi-positive, but not semi-ample. 
We note that $\kappa (X, \mathcal L)=0$ since 
$H^0(X, \mathcal L^{\otimes k})=H^0(C, S^k(\mathcal E))=H^0(C, \mathcal 
O_C)=\mathbb C$ for every $k\geq 0$, 
where $S^k(\mathcal E)$ is the $k$-th symmetric 
product of $\mathcal E$.  
We can easily check that 
$$K_X=\pi^*(K_C\otimes \det \mathcal E)\otimes \mathcal L^{\otimes (-2)}
=\pi^*(K_C\otimes L)\otimes \mathcal L^{\otimes (-2)}, $$  
where $\pi:X\to C$ is the projection. 
Let $m$ be an integer with $m\geq 2$. 
Then 
\begin{align*}
H^i(X, K_X\otimes \mathcal L^{\otimes m})&=H^i(X, \pi^*(K_C\otimes L)
\otimes 
\mathcal L^{\otimes (m-2)})
\\&=\overset{m-1}{\underset{k=1}\bigoplus} 
H^i(C, K_C\otimes L^{\otimes k}).
\end{align*} 
Thus, $h^0(X, K_X\otimes \mathcal L^{\otimes m})=(m-1)(g-1)$, 
$h^1(X, K_X\otimes \mathcal L^{\otimes m})=0$, 
and $h^2(X, K_X\otimes \mathcal L^{\otimes m})=0$ for 
$m\geq 2$. 
So, we obtain no interesting results from injectivity 
theorems. Note that $H^2(X, K_X\otimes \mathcal L^{\otimes m})=0$ for 
$m\geq 1$ also follows from the Kawamata--Viehweg vanishing theorem 
since $\nu(X, \mathcal L)=1$ and $\dim X=2$. 
\end{ex}
\begin{ex}[Cutkosky]\label{ex4}
We use the same notation as in Example \ref{ex3}. 
Let $P\in C$ be a closed point. 
We put $\mathcal F:=\mathcal O_C(P)\oplus L$ and 
$Y:=\mathbb P_C(\mathcal F)$. 
Then it is easy to see that $\mathcal M:=\mathcal O_
{\mathbb P_C(\mathcal F)}(1)$ is big and semi-positive, 
but not semi-ample. 
We note that 
$$\underset{m\geq 0}{\bigoplus} H^0(Y, \mathcal M^{\otimes 
m})$$ is not finitely generated. 
For details, see, for example, \cite[Example 2.3.3]{lazarsfeld1}.  
\end{ex}

The following example shows the difference 
between Enoki's injectivity theorem and Koll\'ar's one. 

\begin{ex}
Let $C$ be a smooth projective curve with the 
genus 
$g(C)=g\geq 1$. 
Let $L\in \Pic ^0(C)$ be non-torsion. 
We put $\mathcal E:=\mathcal O_C\oplus L\oplus L^{-1}$, 
$X:=\mathbb P_C(\mathcal E)$, 
and $\mathcal L:=\mathcal O_{\mathbb P_C(\mathcal E)}(1)$. 
It is obvious that $\mathcal E$ has a smooth hermitian 
metric whose 
curvature is Nakano semi-positive. 
Thus, $\mathcal L$ is semi-positive since $\mathcal L$ 
is a quotient line 
bundle of $\pi^*\mathcal E$, where 
$\pi:X\to C$ is the projection. 
In particular, $\mathcal L$ is nef. 
On the other hand, it is not difficult to see 
that $\mathcal L$ is not semi-ample. 
We have $$K_X=\pi^*(K_C\otimes \det \mathcal E)\otimes \mathcal L^{\otimes (-3)}
=\pi^*K_C\otimes \mathcal L^{\otimes (-3)}. $$ 
We can easily check that 
$$
S^m(\mathcal E)=\underset{\begin{subarray}{c}0\leq a+b\leq m\\ 
a, b\geq 0\end{subarray}}\bigoplus L^{a-b}. 
$$ 
Note that the rank of $S^m(\mathcal E)$ is $\frac{1}{2}(m+2)(m+1)$. 
Let $m$ be an integer with $m\geq 3$. 
Then it is easy to see that 
$$H^i(X, K_X\otimes \mathcal L^{\otimes m})=H^i(X, \pi^*K_C\otimes 
\mathcal L^{\otimes (m-3)})=H^i(C, K_C\otimes 
S^{m-3}(\mathcal E)).$$ 
for all $i$. 
We need the following obvious lemma.  
\begin{lem}
We have $h^0(C, K_C)=g$ and $h^1(C, K_C)=1$. 
Moreover, $h^0(C, K_C\otimes L^{\otimes k})=g-1$ and 
$h^1(C, K_C\otimes L^{\otimes k})=0$ for $k\ne 0$. 
\end{lem} 
Therefore, we obtain $H^3(X, K_X\otimes \mathcal L^{\otimes m})=
H^2(X, K_X\otimes \mathcal L^{\otimes m})=0$, 
$h^1(X, K_X\otimes \mathcal L^{\otimes m})=\llcorner \frac{m-3}{2}\lrcorner +1
=\llcorner \frac{m-1}{2}\lrcorner$, 
and 
\begin{align*}
&h^0(X, K_X\otimes \mathcal L^{\otimes m})\\&=g\llcorner \frac{m-1}{2}
\lrcorner+(g-1)\left(\frac{(m+2)(m+1)}{2}-\llcorner \frac{m-1}{2}\lrcorner\right)
\\&=\llcorner\frac{m-1}{2}\lrcorner +(g-1)\frac{(m+2)(m+1)}{2}. 
\end{align*} 
On the other hand, $h^0(X, \mathcal L^{\otimes k})=h^0(C, S^k(\mathcal E))
=\llcorner \frac{k}{2}\lrcorner+1$ for $k\geq 0$. 
Let $s\in |\mathcal L^{\otimes k}|$ be a non-zero holomorphic section of $\mathcal 
L^{\otimes k}$ for $k\geq 0$. 
Then 
$$
\times s: H^1(X, K_X\otimes \mathcal L^{\otimes m})
\to  
H^1(X, K_X\otimes \mathcal L^{\otimes (m+k)})$$ 
is injective by Enoki's injectivity theorem (cf.~\cite[Theorem 0.2]{enoki}, 
\cite[Corollary 1.4]{f2}).
Note that $h^1(X, K_X\otimes \mathcal L^{\otimes m})
=\llcorner \frac{m-1}{2}\lrcorner$ 
and $h^1(X, K_X\otimes \mathcal L^{\otimes (m+k)})=\llcorner 
\frac{m+k-1}{2}\lrcorner$. 
We have $\kappa (X, \mathcal L)=1$ by the above calculation. 
Since $\mathcal L^2\cdot F=1$, where 
$F$ is a fiber of $\pi:X\to C$, 
$\nu(X, \mathcal L)=2$. 
Thus, the nef line bundle $\mathcal L$ is not abundant. 
So, we think that there are no algebraic proofs for the above injectivity theorem. 
Note that $H^3(X, K_X\otimes \mathcal L^{\otimes m})=
H^2(X, K_X\otimes \mathcal L^{\otimes m})=0$ for $m\geq 1$ follows from the 
Kawamata--Viehweg 
vanishing theorem since $\nu(X, \mathcal L)=2$ and 
$\dim X=3$. 
\end{ex}

The following two examples are famous ones due to 
Mumford and Ramanujam. 

\begin{ex}[Mumford]\label{106}
Let us recall the construction of Mumford's 
example (see \cite[Example 10.6]{ample}). 
We use the same notation as in \cite[Example 10.6]{ample}. 
Let $C$ be a smooth projective curve of genus 
$g\geq 2$ over $\mathbb C$. 
Then there exists a stable vector bundle $E$ of rank two and $\deg E=0$ 
such that its symmetric powers $S^m(E)$ are stable for 
all $m\geq 1$. 
We consider the ruled surface $X:=\mathbb P_C(E)$. 
Let $D$ be the divisor corresponding to 
$\mathcal O_X(1)$. 
Since $E$ is a unitary flat vector bundle, $\mathcal L:=\mathcal 
O_X(D)\simeq \mathcal O_{\mathbb P_C(E)}(1)$ is semi-positive 
by $\pi^*E\to \mathcal L\to 0$, where $\pi:X\to C$ is the 
projection. 
We know that $H^0(X, \mathcal L^{\otimes m})=H^0(C, S^m(E))=0$ 
since $S^m(E)$ is stable 
and $c_1(S^m(E))=0$ for every $m\geq 1$. 
Thus, $\kappa (X, \mathcal L)=-\infty$. 
On the other hand, $\mathcal L\cdot \mathcal L=0$ and 
$\mathcal L\cdot C'>0$ for every curve $C'$ on $X$. 
Then $\nu (X, \mathcal L)=1$. 
\end{ex}

\begin{ex}[Ramanujam]\label{ex58}  
Let us recall the construction of Ramanujam's 
example (see \cite[Example 10.8]{ample}). 
We use the same notation as in \cite[Example 10.8]{ample}. 
Let $X$ be the ruled surface obtained in Example 
\ref{106}. 
We assume that $D$ is the divisor given in \cite[Example 10.6]{ample} 
(see Example \ref{106} above). 
Let $H$ be an effective 
ample divisor on $X$. 
We 
define $\overline {X}:=\mathbb P_X(\mathcal O_X(D-H)
\oplus \mathcal O_X)$, and 
let $\pi:\overline X\to X$ be the projection. 
Let $X_0$ be the section of $\pi$ 
corresponding to $\mathcal O_X(D-H)\oplus 
\mathcal O_X\to \mathcal O_X(D-H)\to 0$ and 
$\overline D:=X_0+\pi^*H$. 
We put $\mathcal M:=\mathcal O_{\overline X}(\overline D)$. 
We write 
$\mathcal O_{\overline X}(1)
=\mathcal O_{\mathbb P_X(\mathcal O_X(D-H)\oplus \mathcal O_X)}(1)$. 
Then $\mathcal O_{\overline X}(1)\simeq \mathcal O_{\overline X}(X_0)$. 
Therefore, 
$\mathcal M\simeq \mathcal O_{\overline X}(1)\otimes \pi^*
\mathcal O_X(H)$ and 
$\pi^*(\mathcal O_X(D)\oplus \mathcal O_X(H))\to 
\mathcal M\to 0$. 
Thus, it is easy to see that $\mathcal M$ is semi-positive, nef and big. 
By the construction, $\mathcal M\cdot C'>0$ for every curve $C'$ on $\overline X$. 
However, $\mathcal M$ is not semi-ample since 
$\mathcal M|_{X_0}\simeq \mathcal O_{X_0}(D)$ does not have 
sections on $X_0$. 
In particular, $\underset{m\geq 0}\bigoplus H^0(\overline X, \mathcal O_{\overline X}
(m\overline D))$ is not finitely generated. 
\end{ex}

The author learned the following construction from the referee. 

\begin{ex}\label{ex-new} 
Let $X=\mathbb P_C(\mathcal E)\to C$ be as in Example \ref{ex1}. 
Let $A$ be a very ample divisor on $C$. 
We take a smooth member $B\in |2A|$ and take 
a double cover $\widetilde C\to C$ by $B\sim 2A$. 
We consider the base change diagram 
$$
\xymatrix{
X\ar[d]&\widetilde X\ar[l]_{\varphi}\ar[d]\\
C& \ar[l]\widetilde C
}
$$ 
and $\mathcal A:=\varphi^*\mathcal O_{\mathbb P_C{(\mathcal E)}}(1)
\simeq \mathcal O_{\widetilde X}(\varphi^*D)$. 
The natural map 
$$\alpha:H^1(\widetilde X, K_{\widetilde X}
\otimes \mathcal O_{\widetilde X}(2\varphi^*D))
\to H^1(\widetilde X, K_{\widetilde X}\otimes \mathcal O_{\widetilde X}(3\varphi^*D))
$$ contains 
$$\beta:H^1(X, K_X\otimes \mathcal O_X(2D))\to 
H^1(X, K_X\otimes \mathcal O_X(3D))$$ 
as a direct summand by the construction. 
Since $\beta$ is zero (see Example \ref{ex1}), 
$\alpha$ is not injective. 
By Enoki's injectivity theorem (cf.~\cite[Theorem 0.2]{enoki}, 
\cite[Corollary 1.4]{f2}), 
$\mathcal A$ is not semi-positive. 
We note that $\widetilde C$ is a smooth 
projective curve of genus $g\geq 2$. 
Then there exists a stable vector bundle $F$ of rank two and $\deg F=0$ such that 
its symmetric powers $S^m(F)$ are stable for 
all $m\geq 1$ (cf.~Example \ref{106}). 
We put $Y=\mathbb P_{\widetilde C}(F)$ and $\mathcal B=\mathcal O_{\mathbb P_{\widetilde C}
(F)}(1)$. 
We take the fiber product 
$$
\xymatrix{\widetilde X \ar[d]& V=\widetilde X\times 
_{\widetilde C}Y\ar[l]_{p_1\ \ \ \ \ }\ar[d]^{p_2}\\ 
\widetilde C &\ar[l] Y
}
$$
and put $\mathcal M=p^*_1\mathcal A\otimes p^*_2\mathcal B$. 
Note that $V$ is a smooth projective variety. 
Then it is easy to see that $\mathcal M\cdot C'>0$ 
for every curve $C'$ on $V$ and that $\kappa (V, \mathcal M)=-\infty$. 
We do not know whether $\mathcal M$ is 
semi-positive or not. 
\end{ex}

We close this paper with a question. 

\begin{ques}
In Example \ref{ex-new}, 
are there any smooth hermitian 
metrics on $\mathcal M$ with semi-positive 
curvature? 
\end{ques}
%%%%%%%%%%%%%%%%%%%%%%%%%%%
\ifx\undefined\bysame
\newcommand{\bysame|{leavemode\hbox to3em{\hrulefill}\,}
\fi

\end{document}